\documentclass{amsart}

\usepackage[all]{xy}
 \usepackage[usenames,dvipsnames]{pstricks}
\usepackage{ebezier}
\usepackage[pagewise]{lineno}
\usepackage[colorlinks]{hyperref}

\usepackage{latexsym}\usepackage{amsfonts}\usepackage{amssymb}\usepackage{amsthm}
\usepackage{amsmath}\usepackage{newlfont}\usepackage{graphicx}\usepackage{doc}
\usepackage{color}\usepackage{multicol}\usepackage{epic}\usepackage{eepic}
\usepackage{ebezier}

\newtheorem{thrm}{Theorem}
\newtheorem{lem}[thrm]{Lemma}
\newtheorem{cor}[thrm]{Corollary}
\newtheorem{defn}[thrm]{Definition}
\newtheorem{rem}[thrm]{Remark}
\newtheorem{prop}[thrm]{Proposition}

\newtheorem{ques}[thrm]{Question}
\theoremstyle{theorem}
\newtheorem{thmA}{Theorem}

\include{epsf}

\DeclareMathOperator{\limarr}{\displaystyle \lim_{\longrightarrow}}

\def \Dj{\mbox{\raise0.3ex\hbox{-}\kern-0.4em D}}

\begin{document}

\title{Spectral numbers and manifolds with boundary}

\author{Jelena Kati\'c, Darko Milinkovi\'c, Jovana Nikoli\'c}

\address{Matemati\v{c}ki fakultet, Studentski trg 16, 11000 Belgrade, Serbia}

\email{jelenak@matf.bg.ac.rs, milinko@matf.bg.ac.rs, jovanadj@matf.bg.ac.rs}

\thanks{This work is partially supported by Ministry of Education and Science  of Republic of Serbia Project \#ON174034.}

\begin{abstract} We consider a smooth submanifold $N$ with a smooth boundary in an ambient closed manifold $M$ and assign a spectral invariant $c(\alpha,H)$ to every singular homological class $\alpha\in H_*(N)$ and a Hamiltonian $H$ defined on the cotangent bundle $T^*M$. We also derive certain properties of spectral numbers, for example we prove that spectral invariants $c_\pm(H,N)$ associated to the whole Floer homology $HF_*(H,N:M)$ of the submanifold $N$, are limits of the decreasing nested family of open sets.
\end{abstract}
\maketitle




\section{Introduction}
\subsection{Floer Homology for submanifold with boundary} Let $N$ be a closed submanifold with a smooth boubdary of an ambient closed manifold $M$. Considered a conormal bundle of $\partial N$, $\nu^*(\partial N)$, defined as
$$\nu^*(\partial N)=\{(\mathbf{q},\mathbf{p})\in T^*M\mid \mathbf{q}\in \partial N, \mathbf{p}|_{T_{\mathbf{q}}\partial N}=0\},$$ which is a Lagrangian submanifold of the cotangent bundle $T^*M$. Let $\nu_-^*(\partial N)$ denote the negative part of :
$$\nu^*_-(\partial N):=\{(\mathbf{q},\mathbf{p})\in\nu^*(\partial N)\mid \mathbf p(\vec{n})\ge 0,\;\mbox{for}\;\vec{n}\in TN\;\mbox{inner normal to}\;N\}$$ and define a {\it negative conormal to $\overline{N}$}, as
$$\nu^*\overline{N}:=\nu^*_-(\partial N)\cup\nu^*(\partial N).$$ The set $\nu^*_-\overline{N}$, called a {\it negative conormal to $\overline{N}$}, is a singular Lagrangian submanifold that allows a smooth approximation by exact Lagrangian submanifolds. If $\Upsilon$ is such an approximation, we can define the Floer homology groups for the pair $(O_M,\Upsilon)$, denoted by $HF_*(O_M,\Upsilon:H,J)$ where $O_M$ is the zero section of $T^*M$. After establishing canonical isomorphisms between Floer homology with two different approximations $\Upsilon_1$ and $\Upsilon_2$ (possibly with a change of almost complex structure $J$), we define the Floer homology for $N$ as a direct limit of the Floer groups $HF_*(O_M,\Upsilon:H,J)$, by taking the canonical isomorphisms for morphisms that define the direct limit. Floer homology obtained in this way is isomorphic to the Morse homology of the set $N$, with appropriately chosen Morse function that satisfies certain conditions near the boundary $\partial N$.

This construction when $N$ is an open subset of $M$ with a smooth boundary is done by Kasturirangan and Oh in~\cite{KO} (see also~\cite{NZ} for the sheaf theoretic point of view). Our construction differs from the one in~\cite{KO} for open subsets by the fact that the tubular neighbourhood of $\partial U$ with respect to $U$ is always trivial, which makes the situation of a submanifold of codimension zero simpler.

\subsection{Spectral numbers} Spectral numbers or spectral invariants were originally defined by Viterbo in~\cite{V}, in the case of a cotangent bundle, in terms of generating functions for Lagrangian submanifolds. In~\cite{O1,O2}, Oh used Viterbo's ideas to define spectral invariants as homologically visible critical values of the action functional
$$a_H(x):=\int_x\theta-\int_0^1H(x(t),t)dt,$$ which is well defined in a cotangent bundle. Here $\theta$ is the canonical Liouville form. We assume that the domain of $a_H$ is the set of smooth paths with the ends on Lagrangian submanifolds $L_0$ and $L_1$. Let $L_0$ be the zero section $O_M$ and
$L_1=\phi^1_H(L_0)$, where $\phi^1_H$ is a time--one--map generated by a Hamiltonian $H$.
The set of critical points $\mathrm{Crit}(a_H)$ of $a_H$ consists of Hamiltonian paths starting at $O_M$ and ending on $\phi^1_H(O_M)$. Let
$HF_*^{\lambda}(O_M,\phi^1_H(O_M):H,J)$ denotes the filtrated homology defined via filtrated Floer complex:
$$CF_*^{\lambda}(O_M,\phi^1_H(O_M):H):=\mathbb Z_2\langle\{x\in\mathrm{Crit}(a_H)\mid a_H(x)<\lambda\}\rangle.$$
These homology groups are well defined since the boundary map preserves the filtration:
$$\partial:CF_*^{\lambda}(O_M,\phi^1_H(O_M):H)\to CF_*^{\lambda}(O_M,\phi^1_H(O_M):H),$$
due to the well defined action functional that decreases along its ``gradient flows".
For a singular homology class $\alpha$ define
$$\sigma(\alpha,H):=
\inf\{{\lambda}\in\mathbb{R}\mid F_H(\alpha)\in\mathrm{Im}(\imath_*^{\lambda})\}$$ where
$$\imath_*^{\lambda}:HF_*^{\lambda}(O_M,\phi^1_H(O_M):H,J)
\to HF_*(O_M,\phi^1_H(O_M):H,J)$$ is the homomorphism induced by inclusion and
$$F_H:H_*(M)\to HF_*(O_M,\phi^1_H(O_M):H,J)
$$ is an isomorphism between singular and Floer homology groups. The construction of spectral invariants is done in~\cite{O1} in case of cotangent bundle, and in~\cite{O2} for cohomology classes.

It turns out that Oh's invariants and those of Viterbo are in fact the same, see~\cite{M1,M2}.

Using spectral numbers, Oh derived the non--degeneracy of Hofer's metric for Lagrangian submanifolds, a result earlier proved by Chekanov~\cite{C} using different methods. Another application to Hofer geometry is given in~\cite{M3,M4} in the characterization of geodesics in Hofer's metric for Lagrangian submanifolds of the cotangent bundle via quasi--autonomous Hamiltonians.

Spectral invariants in cotangent bundles were also studied by Monzner, Vichery and Zapolsky in~\cite{MVZ}.

Spectral numbers also appear in contexts different (more general) from cotangent bundles. In~\cite{L}, Leclercq constructed spectral invariants for Lagrangian Floer theory in the case when $L$ is a closed submanifold of a compact (or convex in infinity) symplectic manifold $P$ and
$\omega|_{\pi_2(P,L)}=0,\quad \mu|_{\pi_2(P,L)}=0$, where $\mu$ is the Maslov index.

Schwarz constructed spectral invariants for contractible periodic orbits when $(P,\omega)$ is a symplectic manifold with $\omega|_{\pi_2(P)}=0$ and  $c_1|_{\pi_2(P)}=0$ (see~\cite{Sc}).

Symplectic invariants were further investigated by Eliashberg and Polterovich~\cite{EP}, Polterovich and Rosen~\cite{PR}, Oh~\cite{O4}, Humili\`ere, Leclercq and Seyfaddini ~\cite{HLS}, Lanzat~\cite{La} and also in~\cite{D},~\cite{M1,M2}.

\subsection{Results of the paper} Let $N$ be a compact submanifold with boundary $\partial N$ of a closed manifold $M$. In Section~\ref{section:approx+floer} we describe how to construct a singular Lagrangian submanifold $\nu^*\overline{N}\subset T^*M$ associated to $N$, as well as a smooth exact Lagrangian approximation $\Upsilon$ of $\nu^*\overline{N}$. The construction of a Floer homology assigned to $N$ as a direct limit of Floer homology of the pairs $(O_M,\Upsilon)$ is also done in Section~\ref{section:approx+floer}.

Next, we construct a PSS-type isomorphism between Floer homology for $N$ and Morse homology of $N$. Here we choose a Morse function to be admissible, i.e.\ to have gradient trajectories that behave well near the boundary $\partial N$ (see the definition at the beginning of the Section~\ref{section:PSS}). We also impose some transversality conditions to a Hamiltonian $H$. More precisely, we have the following theorem.

\begin{thmA}\label{thrmPSS-intro} Let $N$ be a compact submanifold of $M$ with a smooth boundary $\partial N$. Let $f_N$ be an admissible Morse function on $N$ and $H$ as in~(\ref{eq:prop+trans_cond}). There exists a PSS type isomorphism
$$\Phi:HM_*(f_N)\stackrel{\cong}{\longrightarrow}HF_*(H,N:M).$$\qed
\end{thmA}

This theorem is reformulated as Theorem~\ref{thm:PSS} and proven in Section~\ref{section:PSS}.

In Section~\ref{section:invariants} we assign spectral invariant to every non-zero homology class $[\alpha]\in HF_*(N)$. We define it via PSS isomorphism and we show that it is a limit of spectral invariants in Floer homology of approximations (Theorem~\ref{thm:limit}). We also prove the continuity of spectral invariants with respect to Hofer's norm (Theorem~\ref{thm:Cont_invar}).

Next, we introduce three pair-of-pants type products in Morse and Floer theory for the submanifold $N$. More precisely, we prove the following theorem (formulated as Theorem~\ref{thm:prod_open} in Section~\ref{section:invariants}.

\begin{thmA} There exist pair-of-pants type products:
$$\begin{aligned}
&\circ:HF_*(H_1,N:M)\otimes HF_*(H_2,N:M)\to HF_*(H_3,N:M)\\
&\cdot:HM_*(f_1,N)\otimes HM_*(f_2,N)\to HM_*(f_3,N)\\
&\star:HM_*(f,N)\otimes HF_*(H,N:M)\to HF_*(H,N:M)
\end{aligned}$$
that turns Floer homology $HF_*(H,N:M)$ into a $HM_*(f,N)-$module. The above products satisfy:
$$\Phi(\alpha\cdot\beta)=\Phi(\alpha)\circ\Phi(\beta).$$ \qed
\end{thmA}

Further, we prove the triangle inequality:
$$
c^3_N([\alpha\cdot\beta],H_1\sharp H_2)\le c^1_N([\alpha],H_1)+c^2_N([\beta],H_2),
$$
where $j$ in $c_N^j$ emphasizes the corresponding Morse function $f_j$. This is done in Theorem~\ref{thrm:tirangle}

Using Albers' chimney maps~\cite{A}, we compare spectral invariants for submanifold with boundary $N$ and periodic orbits in $T^*M$ (Theorem~\ref{thm:chimney}).

Spectral numbers $c_\pm$ are defined as the infimum (respectively supremum) of all $\lambda\in\mathbb{R}$ for which the inclusion map $\imath^\lambda_*$ is surjective (respectively trivial), see the definitions on the pages~\pageref{c_+} and~\pageref{c_-}. Finally, we prove the continuity of spectral numbers $c_{\pm}$ in the following sense.

\begin{thmA}\label{thrm:limit-intro} Let submanifold $N\subset M$ with boundary $\partial N$ be framed (see the definition in the page~\pageref{framed}).
Let $U_n$ be a decreasing sequence of open subsets with smooth boundaries and $\bigcap_nU_n=N$, where $N$ is a smooth submanifold with a smooth boundary. Then $\lim_{n\to\infty}c_+(H,U_n)=c_+(H,N)$.\qed
\end{thmA}

Theorem~\ref{thrm:limit-intro} is also true for $c_-$ and therefore for $\gamma(H,N):=c_+(H,N)-c_-(H,N)$. As a consequence of Theorem~\ref{thrm:limit-intro} we conclude that
$$c_+(H,\partial N)\le c_+(H,N),$$
see Remark~\ref{rem}.

Theorem~\ref{thrm:limit-intro} is formulated and proven as Theorem~\ref{thrm:cont} in Section~\ref{section:nested}. The main step in the proof is the construction of a map from filtered Floer homology for $N$ to filtered Floer homology for $U_n$, which will change the filtration for an arbitrary small amount.

\section{Approximations and Floer homology} \label{section:approx+floer}

Let $M$ be a closed smooth manifold of dimension $n$ and $N\subset M$ a smooth compact $k$-dimensional submanifold with a smooth boundary (i.e.\ $\partial N\subset N$). Let $\mathrm{Col}(\partial N)$ denote a collar neighbourhood of $\partial N$ in $N$. We can assume that $\mathrm{Col}(\partial N)$ is a subset of a tubular neighbourhood $\mathrm{Tb}_M(\partial N)$ in $M$, which is isomorphic to the normal vector bundle over $\partial N$, therefore $\mathrm{Col}(\partial N)$ induces a line bundle over $\partial N$, denote it by $\mathrm{Tb}_N(\partial N)$. Since $\mathrm{Col}(\partial N)\cong\partial N\times [0,+\infty)$, this bundle is trivial:
\begin{equation}\label{eq:triv}
\mathrm{Tb}_N(\partial N)\cong\partial N\times\mathbb{R}.
\end{equation}
Let us denote by $E$ a $(n-k-1)$-rank vector sub-bundle of $\mathrm{Tb}_M(\partial N)$ such that
$$\mathrm{Tb}_M(\partial N)\cong \mathrm{Tb}_N(\partial N)\oplus E$$
and by
$$\tau:\mathrm{Tb}_N(\partial N)\oplus E\to\mathrm{Tb}_N(\partial N),\quad e_1\oplus e_2\mapsto e_1.$$
Let
$$\pi_1:\mathrm{Tb}_N(\partial N)\to\partial N,\quad\pi_2:E\to\partial N$$ be the vector bundle projections. Let us check that $\mathrm{Tb}_M(\partial N)$ is a $(n-k-1)$-rank vector bundle over $\mathrm{Tb}_N(\partial N)$ with $\tau$ as a projection map. For every $e_1\in\mathrm{Tb}_N(\partial N)$, we have
$$\begin{aligned}\tau^{-1}(e_1)&=\{(\pi_1(e_1),x,e_1)\in\mathrm{Tb}_N(\partial N)\oplus E \mid\pi_2(x)=\pi_1(e_1)\}\\
&\cong\{x\in E\mid\pi_2(x)=\pi_1(e_1)\}=\pi_2^{-1}(\pi_1(e_1))\end{aligned}$$ which is a fibre in the bundle $E$. The transition maps are defined via the transition maps of the bundle $E$.

We wish to define a map
$$\tilde{\tau}:T^*(\mathrm{Tb}_M(\partial N))\to T^*(\mathrm{Tb}_N(\partial N)).$$ Let us choose a connection on $\mathrm{Tb}_M(\partial N)$ which determines the horizontal and vertical space
$$T_p\mathrm{Tb}_M(\partial N)=H_p\oplus V_p$$
at any point $p\in\mathrm{Tb}_M(\partial N)$. For $X\in T\mathrm{Tb}_M(\partial N)$ denote by $X_{\mathrm{hor}}$ the horizontal part of $X$. Let $a\in T^*(\mathrm{Tb}_M(\partial N))$, $Y\in T(\mathrm{Tb}_N(\partial N))$. Define
$$\tilde{\tau}(a)(Y):=a\left(\left(\tau^{-1}Y\right)_{\mathrm{hor}}\right),$$ where $\tau^{-1}$ is any right inverse of $\tau$. Although $\tau^{-1}$ is not uniquely defined, the term $\left(\tau^{-1}Y\right)_{\mathrm{hor}}$ does not depend on the choice of $\tau^{-1}$, but only on the chosen connection.

It follows from~(\ref{eq:triv})
\begin{equation}\label{eq:tubularM}
T^*(\mathrm{Tb}_N(\partial N))\cong T^*(\partial N)\times\underbrace{\left(\mathbb{R}\times\mathbb{R}\right)}_A.
\end{equation}

Let $\nu^*(\partial N)$ and $\nu^*N$ denote the conormal bundles of $\partial N$ and $N$ in $T^*M$. Define
$$\nu^*_-(\partial N):=\{(\mathbf{q},\mathbf{p})\in\nu^*(\partial N)\mid \mathbf p(\vec{n})\ge 0,\;\mbox{for}\;\vec{n}\in TN\;\mbox{inner normal to}\;N\}$$
and
$$\nu^*\overline{N}:=\nu^*_-(\partial N)\cup\nu^*N.$$
Let us the describe the set $\nu^*\overline{N}$ in local coordinates. Let $(q_1,\ldots,q_n,p_1,\ldots,p_n)$ be coordinates in a neighbourhood (in $T^*(\mathrm{Tb}_M(\partial N))$) of a point in $\partial N$ such that
\begin{equation}\label{eq:coordinates}
\begin{aligned}
&\bullet (q_1,\ldots,q_n)\;\mbox{ are coordinates in}\; M\\
&\bullet (q_1,\ldots,q_{k-1,}0,\ldots,0)\;\mbox{ are coordinates in }\;\partial N\\
&\bullet (q_1,\ldots,q_k,0,\ldots,0),\, q_k\le 0\;\mbox{ are coordinates in}\; N.
\end{aligned}
\end{equation}
In such coordinates:
\begin{itemize}
\item the set $\nu^*_-(\partial N)$ is described as
$$\{(q_1,\ldots,q_{k-1},\underbrace{0,\ldots,0}_{n-k+1},\underbrace{0,\ldots,0}_{k-1},p_k,\ldots,p_n)\mid p_k\le 0\}$$
\item the set $\nu^*N$ is described as
$$\{(q_1,\ldots,q_{k},\underbrace{0,\ldots,0}_{n-k},\underbrace{0,\ldots,0}_{k},p_{k+1},\ldots,p_n)\mid q_k\le 0\}$$
\item the set $\nu^*\overline{N}$ is described as
$$\begin{aligned}&\{(q_1,\ldots,q_{k-1},\underbrace{0,\ldots,0}_{n-k+1},\underbrace{0,\ldots,0}_{k-1},p_k,\ldots,p_n)\mid p_k\le 0\}\cup\\
&\{(q_1,\ldots,q_{k},\underbrace{0,\ldots,0}_{n-k},\underbrace{0,\ldots,0}_{k},p_{k+1},\ldots,p_n)\mid q_k\le 0\}.
\end{aligned}$$
\end{itemize}
In these coordinates, the set $A$ from~(\ref{eq:tubularM}), is given by $\{(q_k,p_k)\}$. Denote by $\tilde{\pi}_A:T^*(\partial N)\times A\to A$ the canonical projection and define
$$\pi_A:T^*(\mathrm{Tb}_M(\partial N))\to A,\quad \pi_A:=\tilde{\pi}_A\circ\tilde{\tau}.$$
Let $C$ be the singular curve in $A$ defined as:
$$C:=\{(q_k,0)\mid q_k\ge 0\}\cup\{(0,p_k)\mid p_k\ge 0\}.$$
We have
$$\pi_A(\nu^*\overline{N}\cap T^*(\mathrm{Tb}_M(\partial N)))=C.$$

\begin{figure}
\centering
\includegraphics[height=7cm, scale=1.2]{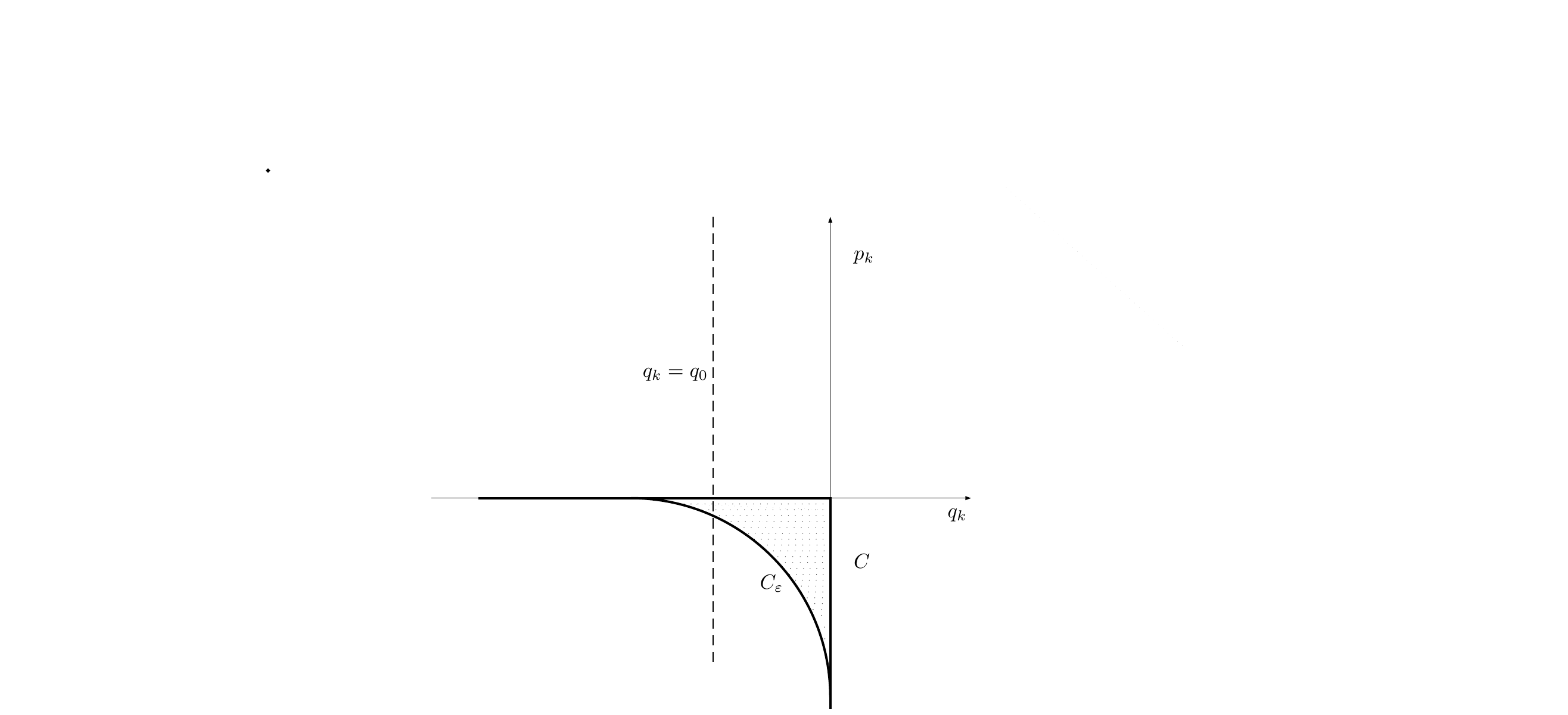}
\centering
\caption{Function $h_{\Upsilon_\varepsilon}$ is the area of the shaded region}
\label{approx2_pic}
\end{figure}

Let $C_1$ be any fixed smooth decreasing curve that coincides with $x$-axis for $x\in(-\infty,-1]$ and with $y$-axis for $y\in(-\infty,-1]$. Let $C_\varepsilon$, for $\varepsilon\in(0,\varepsilon_0)$, be the rescaling of $C_1$ with respect to the cones with vertices in the corners of the curve $C$ (see Figure~\ref{approx2_pic}). We say that the family $\{C_\varepsilon\}_{\varepsilon\in(0,\varepsilon_0)}$ approximates the singular curve $C$ and we say that any curve $C_\varepsilon$ is an approximation of $C$. See also the proof of Theorem 2.3 and Example 2.1 in~\cite{KO}.

Define:
\begin{equation}\label{eq:Ups_eps+coord}
\Upsilon_\varepsilon:=\{(q_1,\ldots,q_{k},\underbrace{0,\ldots,0}_{n-k},\underbrace{0,\ldots,0}_{k},p_k,\ldots,p_n)\mid (q_k,p_k)\in C_\varepsilon\}
\end{equation}

The approximation $\Upsilon_\varepsilon$ is well defined globally. Indeed, let
$$\sigma:T^*(\mathrm{Tb}_N(\partial N))\to T^*(\partial N)$$
denote the projection to the first component in~(\ref{eq:tubularM}). Then it holds
$$\nu^*\overline{N}=\sigma(\nu_-^*(\partial N)\cap\nu^*N)\times C$$ and we can define
\begin{equation}\label{Upsilon_epsilon}
\Upsilon_\varepsilon:=\sigma(\nu_-^*(\partial N)\cap\nu^*N)\times C_\varepsilon.
\end{equation}

Note that $\Upsilon_\varepsilon$ is an exact Lagrangian submanifold of $M$, for every $\varepsilon$. Indeed, if we define the function $h_{\Upsilon_\varepsilon}:\Upsilon_\varepsilon\to\mathbb{R}$ as follows. Let $\pi:T^*M\to M$ denote the canonical projection.
\begin{itemize}\label{def:h_ups}
\item[(a)] on $\nu^*_-\overline{N}\setminus\pi^{-1}(\mathrm{Tb}_M(\partial N))=O_M|_{N\setminus \mathrm{Tb}_M(\partial N)}$:  $h_{\Upsilon_\varepsilon}$ is equal to zero;
\item[(b)] on the intermediate region of $\nu^*_-\overline{N}\cap\pi^{-1}(\mathrm{Tb}_M(\partial N))$: the value of $h_{\Upsilon_\varepsilon}$ at the point
$$(q_1,\ldots,q_{k-1},q_k^0,q_{k+1},\ldots,q_n,p_1,\ldots,p_n)$$ is the negative area of the shaded region in Figure~\ref{approx2_pic} (bounded by $C_\varepsilon$, the $q_k$-axis and the line $q=q_k^0$);
\item[(c)] on $\nu_-^*(\partial N)\cap\Upsilon_\varepsilon$: $h_{\Upsilon_\varepsilon}$ equals the negative area bounded by the $q_k$-axis, the $p_k$-axis and the curve $C_\varepsilon$,
\end{itemize}
then $\theta|_{T\Upsilon_\varepsilon}=dh_{\Upsilon_\varepsilon}$, where $\theta$ is the canonical Liouville form on $T^*M$. Indeed, in the cases $(a)$ and $(c)$ both $\theta|_{T\Upsilon_\varepsilon}$ and $dh_{\Upsilon_\varepsilon}$ equal to zero. Regarding the case $(b)$, if we denote by ${\varphi}_\varepsilon(q_k)$ the function whose graph in $(q_k,p_k)$-plane is the curve $C_\varepsilon$ (see Figure~\ref{approx2_pic}), we have:
$$dh_{\Upsilon_\varepsilon}=\sum_j\frac{\partial h_{\Upsilon_\varepsilon}}{\partial q_j}dq_j+\sum_j\frac{\partial h_{\Upsilon_\varepsilon}}{\partial p_j}dp_j=\varphi(q_k)dq_k,$$ while, on the other hand, we have:
$$\theta|_{T\Upsilon_\varepsilon}=\sum_jp_jdq_j|_{T\Upsilon_\varepsilon}=\varphi(q_k)dq_k$$
(the last equality follows from~(\ref{eq:Ups_eps+coord})).

We will occasionally drop the subscript $\varepsilon$ and denote an approximation only by $\Upsilon$, if the quantity $\varepsilon$ in the approximation is not crucial in that moment.

Floer homology for the pair of exact Lagrangians $(O_M,\Upsilon_\varepsilon)$ is defined in the standard way. Choose a compactly supported Hamiltonian $H:T^*M\times [0,1]\to\mathbb R$ such that
 $$\phi^1_H(O_M)\pitchfork O_M$$ and
 \begin{equation}\label{eq:prop+trans_cond}\phi^1_H(O_M)\cap O_M|_{\partial N}=\emptyset,\quad \phi^1_H(O_M)\pitchfork \nu_-^*\overline{N}.
\end{equation}
Both of the above conditions can be obtained by a generic choice of $H$. The set of the generators of $CF(O_M,\Upsilon_\varepsilon:H)$ consists of Hamiltonian paths
\begin{equation}\label{ham_paths}\dot{x}=X_H(x),\quad x(0)\in O_M,\,x(1)\in \Upsilon_\varepsilon,\end{equation} which are critical points of the {\it effective action functional:}
\begin{equation}\label{eq:ef_act_f-al}
\mathcal{A}_H^{\Upsilon_\varepsilon}(\gamma):=\int \gamma^*\theta-\int_0^1H(\gamma(t),t)dt-h_{\Upsilon_\varepsilon}(\gamma(1)).
\end{equation}

The boundary map $\partial_{J,H}$ is defined by the number of perturbed holomorphic discs with boundary on $O_M$ and $\Upsilon_\varepsilon$:
\begin{equation}\label{d}\left\{\begin{array}{l}
u:\mathbb{R}\times[0,1]\to T^*M\\
\frac{\partial u}{\partial s}+J_\varepsilon\left(\frac{\partial u}{\partial t}-X_H(u)\right)=0\\
u(s,0)\in O_M,\;u(s,1)\in\Upsilon_\varepsilon.
\end{array}\right.
\end{equation}

Floer homology of the manifold with a boundary $N$ is defined as the direct limit of above Floer homologies for the approximations $\Upsilon_\varepsilon$:
\begin{equation}\label{eq:dir_limit}
HF_*(H,N:M):=\displaystyle \lim_{\longrightarrow} HF_*(O_M,\Upsilon_\varepsilon:H,J_\varepsilon),
\end{equation}
after defining an appropriate preorder and a directed set. Take the set of all approximations and define a relation $\le$ on it by
$$\Upsilon_{\varepsilon_1}\le\Upsilon_{\varepsilon_2}\;\Longleftrightarrow\;
\varphi_{\varepsilon_1}\le \varphi_{\varepsilon_2}\quad\mbox{on}\quad U,$$ where
$h_{\Upsilon_{\varepsilon_1}}=\varphi_{\varepsilon_1}\circ\pi$ and $h_{\Upsilon_{\varepsilon_2}}=\varphi_{\varepsilon_2}\circ\pi$, and $\pi:T^*M\to M$ is the canonical projection.
It is not hard to check that the relation $\le$ is a preorder that turns the set of all approximations into a directed set. Since $h_{\Upsilon_{\varepsilon}}$ measures the negative area between $\Upsilon_{\varepsilon}$ and $\nu_-^*\overline{U}$,
$\Upsilon_1\le\Upsilon_2$ means that $\Upsilon_2$ is a better approximation than $\Upsilon_1$.

The connecting morphism that define the direct limit
$$\mathbf F_{12}:HF_*(O_M,\Upsilon_1:H,J_1)\to HF_*(O_M,\Upsilon_2:H,J_2)$$
for $\Upsilon_1\le\Upsilon_2$ is a canonical homomorphism constructed using the standard cobordism arguments (see~\cite{KO} or~(\ref{eq:mon_hom}) below).

\section{Morse homology and PSS isomorphism}\label{section:PSS}

We say that a Morse function $f_N$ defined on $N$ is {\it admissible} if the following holds:
\begin{itemize}
\item[(i)] there are no critical points of $f$ in $\partial N$
\item[(ii)] for all $p\in\partial N$, $df(\vec{n})<0$ where $\vec{n}\in T_pN$ is an inner normal to $N$.
\end{itemize}

We can extend $f_N$ to a Morse function $f$ defined on $N\cup\mathrm{Tb}_N(\partial N)$ such that
\begin{itemize}
\item $\mathrm{Crit}(f)=\mathrm{Crit}(f_N)$
\item no negative gradient trajectories starting at some point in $N$ can leave $N$.
\end{itemize}
See the details in~\cite{S}, subsection 4.2.3.

Choose a Riemannian metric $g$ on $M$ such that the pair $(f_N:g)$ is Morse-Smale. For such $f_N$, Morse homology $HM_*(f_N)$ is isomorphic to the singular homology $H_*(N)$ of $N$.

We can extend $f$ to the normal bundle of $N\cup\mathrm{Tb}_N(\partial N)$ in $M$, and then to the whole $M$ to obtain Morse function $f_M:M\to\mathbb{R}$ with the following properties:\label{page:extension}
\begin{itemize}
\item no negative gradient trajectory $\dot{\gamma}=-\nabla f_M$ leaves $N$
\item there are no critical points of $f_M$ in some neighbourhood of $N$ in $M$
\item Morse complex of $f_N$ is identified with a subcomplex of $f_M$.
\end{itemize} Indeed, let
$$\nu:=\nu(N\cup\mathrm{Tb}_N(\partial N))$$ denotes a tubular neighbourhood of $N\cup\mathrm{Tb}_N(\partial N)$ in $M$. This neighbourhood exists since $N\cup\mathrm{Tb}_N(\partial N)$ is an embedded submanifold of $M$ and it is possesses a $(n-k)$-dimensional vector bundle structure of normal bundle (see Corollary 2.3 in~\cite{K}; if necessary, we can shrink the fibres in $\mathrm{Tb}_N(\partial N)$ to obtain that the closure of $N\cup\mathrm{Tb}_N(\partial N)$ embeds as a closed subset of $M$, which is assumed in~\cite{K}). We will identify the tubular neighbourhood $\nu$ with the corresponding normal bundle.

Let $\langle\cdot,\cdot\rangle_\nu$ be a Riemannian metric with associated quadratic form
$$q(\cdot):=\langle\cdot,\cdot\rangle_\nu.$$ Define
$$f_\nu(p,u_p):=f(p)+q(u_p).$$ The function $f_\nu$ is Morse function defined on $\nu$, critical points of $f_\nu$ are precisely critical points of $f$, of the same Morse index. To define an extension $f_M$ of $f_\nu$, we use the following fact.
\begin{lem}{\bf (Lemma 4.15 in~\cite{S})}
Let $M$ be a smooth manifold, let $A$ be closed and $W$ open subset and $f_W\in C^{\infty}(W,\mathbb{R})$ be a Morse function on the submanifold $W$, such that
$$\mathrm{Crit}(f_W)\subset\mathring{A}\subset A\subset W\subset M.$$ Then there exists a smooth Morse function $f_M$ defined on $M$ which extends $f_W$ meaning that $f_W|_A=f_M|_A$.
\end{lem}
For our needs, we take any closed set $A$ satisfying:
$$N\subset\mathring{A}\subset A\subset\nu$$ and $W:=\nu$.
(See also subsections 4.2.1.\ and 4.2.3.\ in~\cite{S} for more details.)

For two Riemannian metrics $g_a$ and $g_b$ there is a canonical isomorphism
$$\mathbf G_{ab}:HM_*(f_N:g_a)\to HM_*(f_N:g_b)$$
satisfying
$$\mathbf G_{ac}=\mathbf G_{bc}\circ \mathbf G_{ab},\quad\mathbf  G_{aa}=\mathrm{Id}.$$ This functoriality allows us to consider the set $\{HM_*(f_N:g_a)\}$ as a directed system, indexed by the directed set of all generic Riemannian metrics $\{g_a\}$ with the full relation as a preorder. We define Morse homology $HM_*(f_N)$ as a direct limit:
$$HM_*(f_N):=\displaystyle \lim_{\longrightarrow} HM_*(f_N:g_a):=
\bigsqcup_sHM_*(f_N:g_s)/\sim$$
where
$$p_a\sim p_b\quad\Leftrightarrow\quad\mathbf{G}_{ac}(p_a)=\mathbf{G}_{bc}(p_b)$$ for some $c$. The set $HM_*(f_N)$ obviously has a vector space structure and is isomorphic to all $HM_*(f_N:g_s)$.

Concerning the isomorphism between Floer homology $HF_*(H,N:M)$ and Morse homology $HM_*(f_N)$, we have the following theorem.

\begin{thrm}\label{thm:PSS} Let $f_N$ be admissible and $H$ as in~(\ref{eq:prop+trans_cond}). There exists a PSS type isomorphism
$$\Phi:HM_*(f_N)\stackrel{\cong}{\longrightarrow}HF_*(H,N:M).$$
\end{thrm}

\noindent{\it Proof.} We first define the PSS homomorphism for the approximations. Let us fix the approximation $\Upsilon$. Let $p\in N$ be a critical point of $f_N$ and $x$ be a Hamiltonian path with $x(0)\in O_M$, $x(1)\in\Upsilon$.
The PSS homomorphism is defined via the number of mixed objects
$$\begin{aligned}&\mathcal{M}(p,x):=
\mathcal{M}(p,x,O_M,\Upsilon:f,H,J,g):=\\
&\left\{(\gamma,u)\left|\begin{array}{l}
\gamma:(-\infty,0]\to N,\; u:[0,+\infty)\times[0,1]\to T^*M\\
\dot{\gamma}(s)=-\nabla_g f(\gamma(s)) \\
\frac{\partial u}{\partial s}+
J(\frac{\partial u}{\partial t}-X_{\rho_RH}(u))= 0\\
u(s,0), u(0,t)\in O_M,\,u(s,1)\in\Upsilon  \\
\gamma(-\infty)=p,\,
u(+\infty,t)=x(t) \\
u(0,1)=\gamma(0)\\
\end{array}\right.\right\}\end{aligned}$$ where $\rho_R:[0,+\infty)\rightarrow{\mathbb R}$ is a smooth function
such that
$$\rho_R(s)=\begin{cases} 1, & s\ge R \\ 0, & s\le R-1. \end{cases}
$$

The loss of the compactness of the manifolds $\mathcal{M}(p,x)$ manifests in a breaking of gradient trajectories or holomorphic strips. The appearance of a bubbling is excluded due to the exact Lagrangian boundary conditions. This in particular means that the zero dimensional manifolds $\mathcal{M}(p,x)$ are compact, so the map
$$p\mapsto\sum_{x\in CF_*(O_M,\Upsilon:H)}n(p,x)x,\quad n(p,x):=\,\sharp\,\mathcal {M}(p,x)\pmod 2,$$ is well defined on the chain level. On the other hand, the description of the boundary of one-dimensional manifolds of the same type provides the fact that the above map descends on the homology level:
$$\Phi^{\Upsilon}:HM_*(f_N:g)\to HF_*(O_M,\Upsilon:H,J).$$
The map $\Phi^\Upsilon$ turns out to be an isomorphism; by standard cobordism arguments one can prove that its inverse is defined via the number of mixed objects of the type $(u,\gamma)$. All the above holds for generic choices of almost complex structure $J$ and Riemannian metric $g$.
The details of the above construction are standard (see~\cite{KMN} for the case of open subset of the base).

In order to prove that $\Phi^\Upsilon$ defines a map between Morse and Floer homologies defined as the direct limits of approximations, one needs to check if it commutes with the morphisms defining the direct limits. This means that, if we denote:
$$\begin{aligned}\Phi^{a}:HM_*(f_N:g_a)\to HF_*(O_M,\Upsilon_a:H,J_a),\\\Phi^{b}:HM_*(f_N:g_b)\to HF_*(O_M,\Upsilon_b:H,J_b),
\end{aligned}
$$ for two approximations $\Upsilon_a$ and $\Upsilon_b$ and (generically chosen) almost complex structures $J_a$ and $J_b$ and Riemannian metrics $g_a$ and $g_b$ it holds
$$\Phi^b\circ\mathbf{G}_{ab}=\mathbf{F}_{ab}\circ\Phi^a.$$
 This also can be done using suitable one-dimensional auxiliary manifolds, similarly to~\cite{KMN}. Let us denote by
 $$\Phi:HM_*(f_N)\to HF_*(H,N:M)$$ the induced homomorphism on the direct limits.

Now if we denote by $\Psi^a:=(\Phi^a)^{-1}$ we have the same arguments as above for the maps $\Psi^a$. Therefore we have a well defined homomorphism
$$\Psi:HF_*(H,N:M)\to HM_*(f_N).$$ From
$$\Psi^a\circ\Phi^a=\mathrm{Id}_{HM_*(f_N:g_a)},\quad\Phi^a\circ\Psi^a=\mathrm{Id}_{HF_*(O_M,\Upsilon_a:H,J_a)}$$ we deduce:
$$\Psi\circ\Phi=\mathrm{Id}_{HM_*(f_N)},\quad\Phi\circ\Psi=\mathrm{Id}_{HF_*(H,N:M)}.$$ \qed

\section{Spectral invariants}\label{section:invariants}

In order to define spectral invariants, we need to establish filtered Floer homology, both for approximations and for direct limit homologies.

Recall that the filtered Floer homology groups for approximations are defined as homology groups of the filtered chain complex
$$CF^\lambda_*(O_M,\Upsilon:H):=\{x\in CF_*(O_M,\Upsilon:H)\mid \mathcal{A}^\Upsilon_H(x)<\lambda\},$$ where the effective action functional $\mathcal{A}^\Upsilon_H$ is defined in~(\ref{eq:ef_act_f-al}).
Since $\mathcal{A}^\Upsilon_H$ decreases along the strips that define the boundary operator
$$\partial_{J,H}:CF_*(O_M,\Upsilon:H)\to CF_{k-1}(O_M,\Upsilon:H),$$
the boundary operator descends to $CF^{\lambda}_*(O_M,\Upsilon:H)$ and defines
$$\partial_{J,H}^\lambda:CF_*^\lambda(O_M,\Upsilon:H)\to CF^\lambda_*(O_M,\Upsilon:H).$$
We denote the corresponding homology groups by $HF_*^\lambda(O_M,\Upsilon:H,J_\Upsilon)$ and by
$$\imath^{\lambda}_{\Upsilon*}:HF_*^\lambda(O_M,\Upsilon:H,J_\Upsilon)\to HF_*(O_M,\Upsilon:H,J_\Upsilon)$$ the homomorphism
induced by the inclusion map $\imath_\Upsilon^{\lambda}$. As usually, for $\alpha\in HM_*(f_N:g_\Upsilon)\setminus\{0\}$ define
$$c_\Upsilon(\alpha,H):=\inf\{\lambda\mid\Phi^\Upsilon(\alpha)\in\mathrm{Im}(\imath^\lambda_{\Upsilon*})\}.$$

The proof of the following proposition is the same as the proof of Lemma 2.6 in~\cite{MVZ}.
\begin{prop}\label{prop:idep+H} If $H$ and $K$ are two Hamiltonian functions satisfying~(\ref{eq:prop+trans_cond}) such that $\phi_H^1=\phi_K^1$, then it holds $c_\Upsilon(\alpha,H)=c_\Upsilon(\alpha,K)$.
\end{prop}

The spectral invariants for approximations actually do not depend on the given approximation, assuming that it is good enough. More precisely:

\begin{prop}\label{prop:limit} Let $\alpha\in HM_*(f_N:g_\Upsilon)\setminus\{0\}$. (Here $g_\Upsilon$ is a generic metric chosen such that the corresponding PSS for approximation $\Upsilon$ is well defined.)
There exists an approximation $\widetilde{\Upsilon}$ such that all invariants $c_{\widetilde{\Upsilon}}(\mathbf{G}_{\Upsilon\widetilde{\Upsilon}}(\alpha),H)$ are equal for all $\widetilde{\Upsilon}$ with $\widetilde{\Upsilon}\le\overline{\Upsilon}$. Here $\mathbf{G}_{\Upsilon\widetilde{\Upsilon}}$ is a canonical isomorphism between $HM_*(f_N:g_\Upsilon)$ and $HM_*(f_N:g_{\widetilde{\Upsilon}})$.
\end{prop}

Proposition~\ref{prop:limit} will follow directly from Theorem~\ref{thm:limit} that we prove below.

\medskip

Therefore, one natural way to define spectral invariants for Floer homology of a singular Lagrangian $\nu^*\overline{N}$ is to take a limit of spectral invariant for approximations. Let us also mention another way to define them, via PSS isomorphism.

The filtered Floer homology for a manifold with boundary $N$ is again defined as a direct limit of corresponding filtered homologies of approximations. The direct limit homomorphisms $\mathbf{F}_{ab}$ are defined via the number of elements in
\begin{equation}\label{eq:mon_hom}
\mathcal{M}(x,y:\widetilde{\Upsilon}_s):=
\left\{ u\left|
\begin{array}{l}
u:\mathbb{R}\times[0,1]\rightarrow T^*M\\
\frac{\partial u}{\partial s}+
\widetilde J_s(\frac{\partial u}{\partial t}-X_H(u))= 0 \\
u(s,0)\in O_M,\,u(s,1)\in\widetilde{\Upsilon}_s\\
u(-\infty,t)=x(t),\;u(+\infty,t)=y(t).
\end{array}\right.\right\}
\end{equation}
Here $\widetilde\Upsilon_s$ is a monotone homotopy for $s\in\mathbb{R}$ such that
$$
\widetilde\Upsilon_s=\begin{cases}
\Upsilon_a, &s\le-R\\
\Upsilon_b, &s\ge R.
\end{cases}$$
(By monotone homotopy we mean $s_1\le s_2\Rightarrow\Upsilon_{s_1}\le\Upsilon_{s_2}$.) The corresponding action functional $\mathcal{A}_H^{\widetilde{\Upsilon}_s}$ decreases along perturbed holomorphic strips that define $\mathbf{F}_{ab}$, therefore
$$\mathcal{A}_H^{\Upsilon_b}(u(+\infty,t))\le \mathcal{A}_H^{\Upsilon_a}(u(-\infty,t))$$ whenever there exists an $u\in\mathcal{M}(x,y:\widetilde{\Upsilon}_s)$. In particular, the homomorphisms $\mathbf{F}_{ab}$ descend to the filtered chain complex. By standard arguments one shows that they are also well defined on filtered homology groups:
$$\mathbf{F}_{ab}^\lambda:HF_*^{\lambda}(O_M,\Upsilon_a:H,J_a)\to HF_*^{\lambda}(O_M,\Upsilon_b:H,J_b).$$
Now we define the filtered Floer homology for $N$ as the direct limit:
$$HF_*^{\lambda}(H,N:M):=\limarr HF_*^{\lambda}(O_M,\Upsilon_s:H,J_s).$$ It is easy to see that
$$\mathbf{F}_{ab}\circ\imath^\lambda_{\Upsilon_a*}=\imath^\lambda_{\Upsilon_b*}\circ\mathbf{F}_{ab}^\lambda,$$
where
$$\imath^\lambda_{\Upsilon_a*}:HF_*^{\lambda}(O_M,\Upsilon_a:H,J_a)\to HF_*(O_M,\Upsilon_a:H,J_a)$$ denotes the inclusion-induced map for the approximations. Thus the induced inclusion maps
$$\imath^\lambda_*:HF_*^{\lambda}(H,N:M)\to HF_*(H,N:M)$$ are also well defined.

\begin{defn} Let $[\alpha]\in HM_*(f_N)\setminus\{0\}$. By class in $[\alpha]$ we assume the class in the direct limit, so $\alpha\in HM_*(f_N:g)$, for some metric $g$. A {\it spectral invariant for a manifold $N$ with boundary} is defined as
\begin{equation}\label{eq:def_inv_c_open}
c_N([\alpha],H):=\inf\{\lambda\mid \Phi([\alpha])\in\mathrm{Im}(\imath_*^\lambda)\}.
\end{equation}
\end{defn}

A natural question is weather $c_N([\alpha],H)$ equals to the limit of $c_\Upsilon(\alpha,H)$ in the sense of Proposition~\ref{prop:limit}. The answer is affirmative and it is given in the following theorem.

\begin{thrm}\label{thm:limit} Let $\alpha\in HM_*(f_N:g_\Upsilon)\setminus\{0\}$. There exists an approximation $\widetilde{\Upsilon}$ such that, for all $\overline{\Upsilon}$ with $\widetilde{\Upsilon}\le\overline{\Upsilon}$ it holds:
$$c_{\overline{\Upsilon}}(\mathbf{G}_{\Upsilon\overline{\Upsilon}}(\alpha),H)=c_N([\alpha],H).$$
\end{thrm}

\noindent{\it Proof.} Let $[\alpha]\in HM_*(f_N)\setminus\{0\}$, $\lambda\in{\mathbb R}$ and $[x]\in HF_*^\lambda(H,N:M)$ be such that
$$\Phi([\alpha])=\imath^\lambda_*([x]).$$
Let $\alpha\in HM_*(f_N:g_\Upsilon)$, $x\in HF_*^\lambda(O_M,\Upsilon':H,J_{\Upsilon'})$
for some $\Upsilon$ and $\Upsilon'$. Since
$$
\Phi([\alpha])=[\Phi^\Upsilon(\alpha)]=[\imath^\lambda_{\Upsilon'*}(x)]=\imath^\lambda_*[x],
$$
we have $\mathbf F_{\Upsilon\widehat{\Upsilon}}(\Phi^\Upsilon(\alpha))=\mathbf F_{\Upsilon'\widehat{\Upsilon}}(\imath^\lambda_{\Upsilon'*}(x))$. We have the following commutative diagram:
\begin{equation}\label{eq:2diagrams}
\begin{array}{lllllllll}
\cdots&\rightarrow&HM_*(f_N:g_{\Upsilon})&\stackrel{\mathbf G_{\Upsilon\widetilde{\Upsilon}}}{\rightarrow}&HM_*(f_N:g_{\widetilde\Upsilon})&\stackrel{\mathbf G_{\widetilde{\Upsilon}\widehat\Upsilon}}{\rightarrow}&HM_*(f_N:g_{\widehat\Upsilon})&\rightarrow&\cdots\\
&&\downarrow \Phi^{\Upsilon} &&\downarrow \Phi^{\widetilde\Upsilon}&&\downarrow \Phi^{\widehat\Upsilon}&&\\
\cdots&\rightarrow&HF_*(\Upsilon:J_\Upsilon)&\stackrel{\mathbf F_{\Upsilon\widetilde{\Upsilon}}}{\rightarrow}&HF_*(\widetilde\Upsilon:J_{\widetilde{\Upsilon}})&\stackrel{\mathbf F_{\widetilde{\Upsilon}\widehat\Upsilon}}{\rightarrow}&HF_*(\widehat\Upsilon:J_{\widehat{\Upsilon}})&\rightarrow&\cdots\\
&&\uparrow \imath^{\lambda}_{\Upsilon*} &&\uparrow \imath^{\lambda}_{\widetilde\Upsilon*}&&\uparrow \imath^{\lambda}_{\widehat\Upsilon*}&&\\
\cdots&\rightarrow&HF_*^\lambda(\Upsilon:J_\Upsilon)&\stackrel{\mathbf F_{\Upsilon\widetilde{\Upsilon}}^\lambda}{\rightarrow}&HF_*^\lambda(\widetilde\Upsilon:J_{\widetilde{\Upsilon}})&\stackrel{\mathbf F_{\widetilde{\Upsilon}\widehat\Upsilon}^\lambda}{\rightarrow}&HF_*^\lambda(\widehat\Upsilon:J_{\widehat{\Upsilon}})&\rightarrow&\cdots
\end{array}
\end{equation}
 In~(\ref{eq:2diagrams}) we abbreviated $HF_*(O_M,\Upsilon:H,J_\Upsilon)$ to $HF_*(\Upsilon:J_\Upsilon)$ and so on. From this commutativity we read
$$
\Phi^{\widehat{\Upsilon}}(\mathbf G_{\Upsilon\widehat{\Upsilon}}(\alpha))=\mathbf F_{\Upsilon\widehat{\Upsilon}}(\Phi^\Upsilon(\alpha))=\mathbf F_{\Upsilon'\widehat{\Upsilon}}(\imath^\lambda_{\Upsilon'*}(x))
=\imath^{\lambda}_{\widehat{\Upsilon}*}(\mathbf F^\lambda_{\Upsilon'\widehat{\Upsilon}}(x)),
$$ implying $\Phi^{\widehat{\Upsilon}}(\mathbf G_{\Upsilon\widehat{\Upsilon}}(\alpha))\in\mathrm{Im}(\imath^\lambda_{\widehat{\Upsilon}*})$.
We conclude
\begin{equation}\label{eq:spect1}
c_{\widehat{\Upsilon}}(\mathbf G_{\Upsilon\widehat{\Upsilon}}(\alpha),H)\leq c_N([\alpha],H).
\end{equation}
If we take $\alpha\in HM_*(f_N:g_{\Upsilon})\setminus\{0\}$ and $\lambda\in{\mathbb R}$ such that
$$\Phi^\Upsilon(\alpha)\in\mathrm{Im}(\imath^\lambda_{\Upsilon*}),$$
then $\Phi^\Upsilon(\alpha)=\imath^\lambda_{\Upsilon*}(x)$ for some $x\in HF_k^\lambda(\Upsilon:H,J_\Upsilon)$. Therefore, we have
$$
\Phi([\alpha])=[\Phi^\Upsilon(\alpha)]=[\imath^\lambda_{\Upsilon*}(x)]=\imath^\lambda_*[x],
$$
so we obtain the inequality
\begin{equation}\label{eq:spect2}
c_N([\alpha],H)\leq c_\Upsilon(\alpha,H).
\end{equation}
The elements $\alpha$ and $\mathbf G_{\Upsilon\widehat{\Upsilon}}(\alpha)$ represent the same element in the quotient space $HM_*(f_N)$. From (\ref{eq:spect1}) and (\ref{eq:spect2}) we have
\begin{equation}\label{eq:ineqaux}
c_{\widehat{\Upsilon}}(\mathbf G_{\Upsilon\widehat{\Upsilon}}(\alpha),H)\leq c_N([\alpha],H)=c_N([\mathbf G_{\Upsilon\widehat{\Upsilon}}(\alpha)],H)\leq c_{\widehat{\Upsilon}}(\mathbf G_{\Upsilon\widehat{\Upsilon}}(\alpha),H),
\end{equation}
so all inequalities become equalities.
It also holds:
\begin{equation}\label{eq:inequ_inv}
c_{\widetilde{\Upsilon}}(\mathbf G_{\Upsilon\widetilde{\Upsilon}}(\alpha),H)\leq c_{\Upsilon}(\alpha,H),
\end{equation}
for every ${\Upsilon}\le\widetilde{\Upsilon}$. Indeed, if $\Phi^\Upsilon(\alpha)\in\mathrm{Im}(\imath^\lambda_{\Upsilon*})$, for $\alpha\neq 0\in HM_*(f_N;g_\Upsilon)$, we have
\begin{equation}\label{eq:auxili}
\Phi^\Upsilon(\alpha)=\imath^\lambda_{\Upsilon*}(x)\;\Rightarrow \;\mathbf F_{\Upsilon\widetilde{\Upsilon}}(\Phi^\Upsilon(\alpha))=\mathbf F_{\Upsilon\widetilde{\Upsilon}}(\imath^\lambda_{\Upsilon*}(x)),
\end{equation} so, from the commutativity~(\ref{eq:2diagrams}) it follows
$$\imath^\lambda_{\widetilde\Upsilon*}(\mathbf F_{\Upsilon\widetilde{\Upsilon}}(x))=\mathbf F_{\Upsilon\widetilde{\Upsilon}}(\imath^\lambda_{\Upsilon*}(x))\stackrel{(\ref{eq:auxili})}{=}\mathbf F_{\Upsilon\widetilde{\Upsilon}}(\Phi^\Upsilon(\alpha))=\Phi^{\widetilde\Upsilon}(\mathbf G_{\Upsilon\widetilde{\Upsilon}}(\alpha)).
$$
This means that
$$\Phi^{\widetilde\Upsilon}(\mathbf G_{\Upsilon\widetilde{\Upsilon}}(\alpha))\in\mathrm{Im}(\imath^\lambda_{\widetilde\Upsilon*}),$$ so~(\ref{eq:inequ_inv}) holds.

From~(\ref{eq:ineqaux}) and~(\ref{eq:inequ_inv}) we conclude that $c_{\widehat{\Upsilon}}(\mathbf{G}_{\Upsilon\hat\Upsilon}(\alpha),H)$ becomes equal to $c_N ([\alpha], H)$, for every $\widehat{\Upsilon}$ with $\widetilde\Upsilon\le\widehat\Upsilon$.
\qed

Proposition~\ref{prop:limit} now follows directly from Theorem~\ref{thm:limit}.

An immediate consequence of Proposition~\ref{prop:idep+H} and Theorem~\ref{thm:limit} is the following.
\begin{cor} If $H$ and $G$ satisfy~(\ref{eq:prop+trans_cond}) and $\phi_H^1=\phi_G^1$, then $c_N([\alpha],H)=c_N([\alpha],G)$.
\end{cor}

\subsection{Continuity}

Spectral invariants are continuous with respect to the Hofer norm.

\begin{thrm}\label{thm:Cont_invar} Let $\|\cdot\|$ denotes the Hofer's norm:
$$\|H\|:=\int_0^1[\max_xH(x,t)-\min_xH(x,t)]dt.$$
Relative spectral invariants for $N$
$$C_N([\alpha],H):=c_N([\alpha],H)-c_N(1,H)
$$
are continuous with respect to $\|\cdot\|:$
$$
|C_N([\alpha],H)-C_N([\alpha],H')|\le\|H-H'\|.
$$
Here $1$ denotes the generator of zero homology group $HM_0(f_N)$.
\end{thrm}
\noindent{\it Proof.} The first step is to prove the continuity of relative invariants for approximations:
$$C_\Upsilon(\alpha,H):=c_\Upsilon(\alpha,H)-c_\Upsilon(1,H).$$ The proof of theorem then follows from the above inequality and Theorem~\ref{thm:limit}.
To prove the continuity result for the relative invariants for approximation, we consider the linear homotopy $$H^s=(1-s)H+sH'=H+s(H'-H).$$ The canonical isomorphism $S^{\Upsilon}_{H,H'}$ is defined by a number of the holomorphic strips that connect a generator $x$ of $CF_*(O_M,\Upsilon:H,J_\Upsilon)$ and a generator $y$ of $CF_*(O_M,\Upsilon:H',J_\Upsilon)$:
$$\begin{aligned}&\mathcal{M}(x,y,O_M,\Upsilon:H,H',J_\Upsilon):=\\
&\left\{u:{\mathbb R}\times[0,1]\to T^*M\left|\begin{array}{l}
\frac{\partial u}{\partial s}+
J_\Upsilon(\frac{\partial u}{\partial t}-X_{H^s}(u))= 0\\
u(s,0)\in O_M,\,u(s,1)\in\Upsilon,  \\
u(-\infty,t)=x(t),\, u(+\infty,t)=y(t)\\
\end{array}\right.\right\}.\end{aligned}$$
If there exists $u\in\mathcal{M}(x,y,O_M,\Upsilon:H,H',J_\Upsilon)$, for the linear homotopy $H^s$, then by direct computation we see that it holds
\begin{equation}\label{eq:aux_est}
\mathcal{A}_{H'}^\Upsilon(y)-\mathcal{A}_{H}^\Upsilon(x)=\int_{-\infty}^{+\infty}\frac{d}{ds}\mathcal{A}_{H^s}^\Upsilon(u(s,\cdot))\leq E_+(H-H'):=\int_0^1\max_x(H-H')\,dt.
\end{equation}
If the linear homotopy is not regular, we can approximate it by a $C^1$-close regular homotopy $H^s=H+\sigma(s)(H'-H)$, and obtain:
$$\mathcal{A}_{H'}^\Upsilon(y)-\mathcal{A}_{H}^\Upsilon(x)\leq E_+(H-H')+\varepsilon$$ for any $\varepsilon>0$, so the estimate~(\ref{eq:aux_est}) holds for a regular homotopy $H^s$.
We have
\begin{equation}\label{ineq_inv}S_{H,H'}^\Upsilon(x)=\sum\pm y_j
\Rightarrow\mathcal{A}_{H'}^\Upsilon(y_j)
\leq\mathcal{A}_{H}^\Upsilon(x)+E_+(H-H').
\end{equation}
For $x\in HF_*(O_M,\Upsilon:H,J_\Upsilon)$ define
$$
\sigma_\Upsilon(x,H):=\inf\{\lambda\in{\mathbb R}\,|\,x\in \imath_{\Upsilon,H*}^\lambda)\}.
$$
Obviously, it holds:
$$c_\Upsilon(\alpha,H)=\sigma_\Upsilon(\Phi^\Upsilon_H(\alpha),H).$$
It follows from~(\ref{ineq_inv}):
$$
\sigma_\Upsilon(S_{H,H'}^\Upsilon(x),H')\leq\sigma_\Upsilon(x,H)+E_+(H-H').
$$
From $S_{H,H'}^\Upsilon\circ\Phi^\Upsilon_H=\Phi^\Upsilon_{H'}$ we conclude
\begin{equation}\label{eq:ineq-alpha}\begin{aligned}
c_\Upsilon(\alpha,H')&=\sigma_\Upsilon(\Phi^\Upsilon_{H'}(\alpha),H')=\sigma_\Upsilon(S_{H,H'}^\Upsilon\circ \Phi^\Upsilon_H(\alpha),H')\\
&\leq\sigma_\Upsilon(\Phi^\Upsilon_H(\alpha),H)+E_+(H-H')\\
&=c_\Upsilon(\alpha,H)+E_+(H-H'),
\end{aligned}\end{equation}
for all $0\neq\alpha\in HM_*(f_N:g_\Upsilon)$. The proof now easily follows from the inequality~(\ref{eq:ineq-alpha}) and the same one applied to $\alpha=1\in HM_0(f_N:g_\Upsilon)$.\qed

\subsection{Products in Morse and Floer homology}

Products in Morse and Floer theory were studied by various authors: Abbondandolo and Schwarz~\cite{AS}, Auroux~\cite{Aur}, Oh~\cite{O2} and also in~\cite{KMS,KMN}.

Here we establish three products of pair of pant type.
\begin{thrm}\label{thm:prod_open}
Let $f_N$ and $H,H_j$, for $j=1,2,3$ be as in Theorem~\ref{thm:PSS}. There exist pair-of-pants type products:
$$\begin{aligned}
&\circ:HF_*(H_1,N:M)\otimes HF_*(H_2,N:M)\to HF_*(H_3,N:M)\\
&\cdot:HM_*(f_1,N)\otimes HM_*(f_2,N)\to HM_*(f_3,N)\\
&\star:HM_*(f,N)\otimes HF_*(H,N:M)\to HF_*(H,N:M)
\end{aligned}$$
that turns Floer homology $HF_*(H,N:M)$ into a $HM_*(f,N)-$module. The above products satisfy:
$$\Phi(\alpha\cdot\beta)=\Phi(\alpha)\circ\Phi(\beta),$$ where $\Phi$ is a PSS isomorphism from Theorem~\ref{thm:PSS}.
\end{thrm}

\begin{figure}
\centering
\includegraphics[width=8cm,height=3cm]{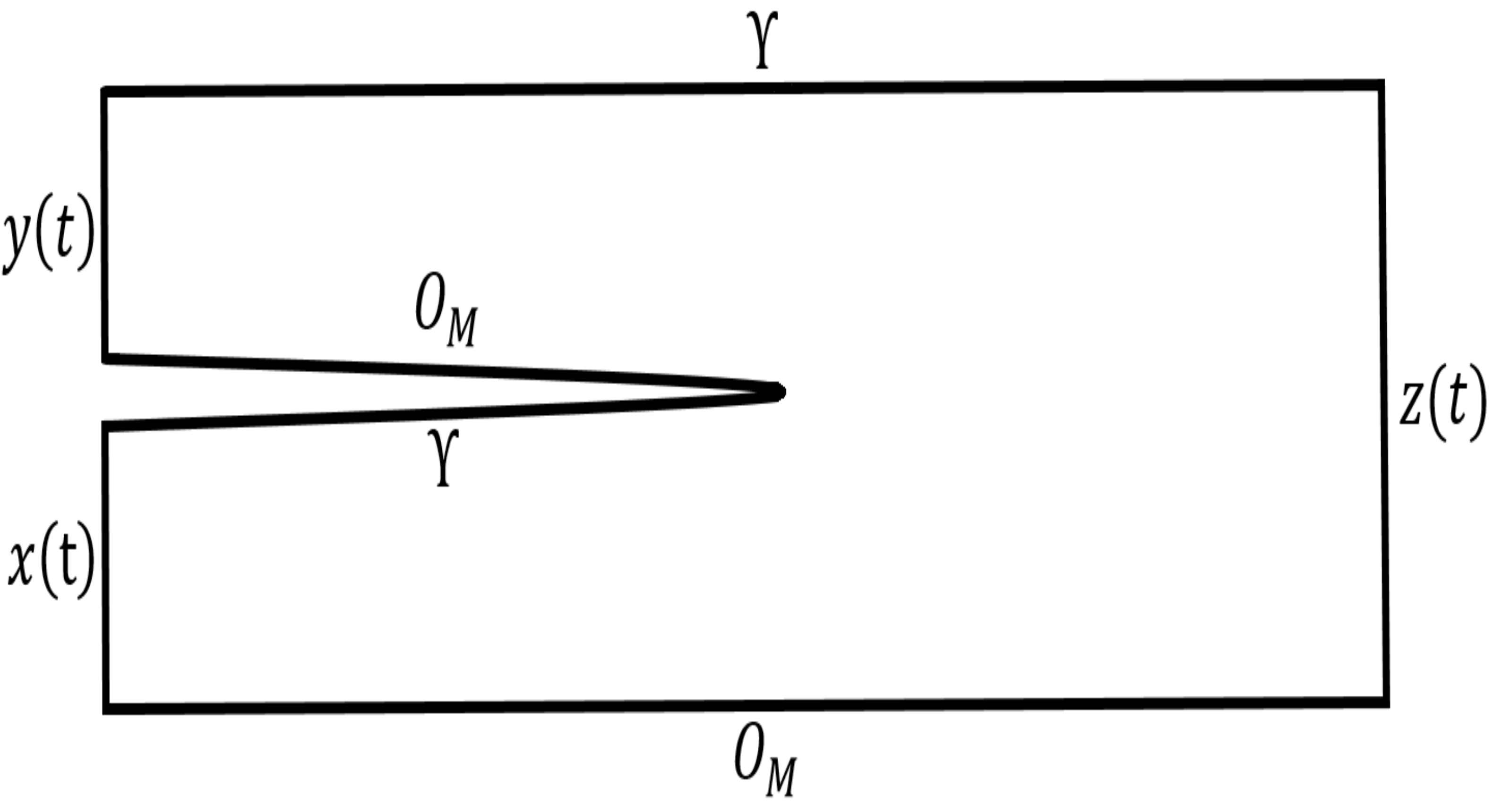}
\centering
\caption{Riemannian surface $\Sigma$ that defines the product $\circ$}
\label{fig:pants}
\end{figure}

\begin{figure}
\centering
\includegraphics[width=6cm,height=2.5cm]{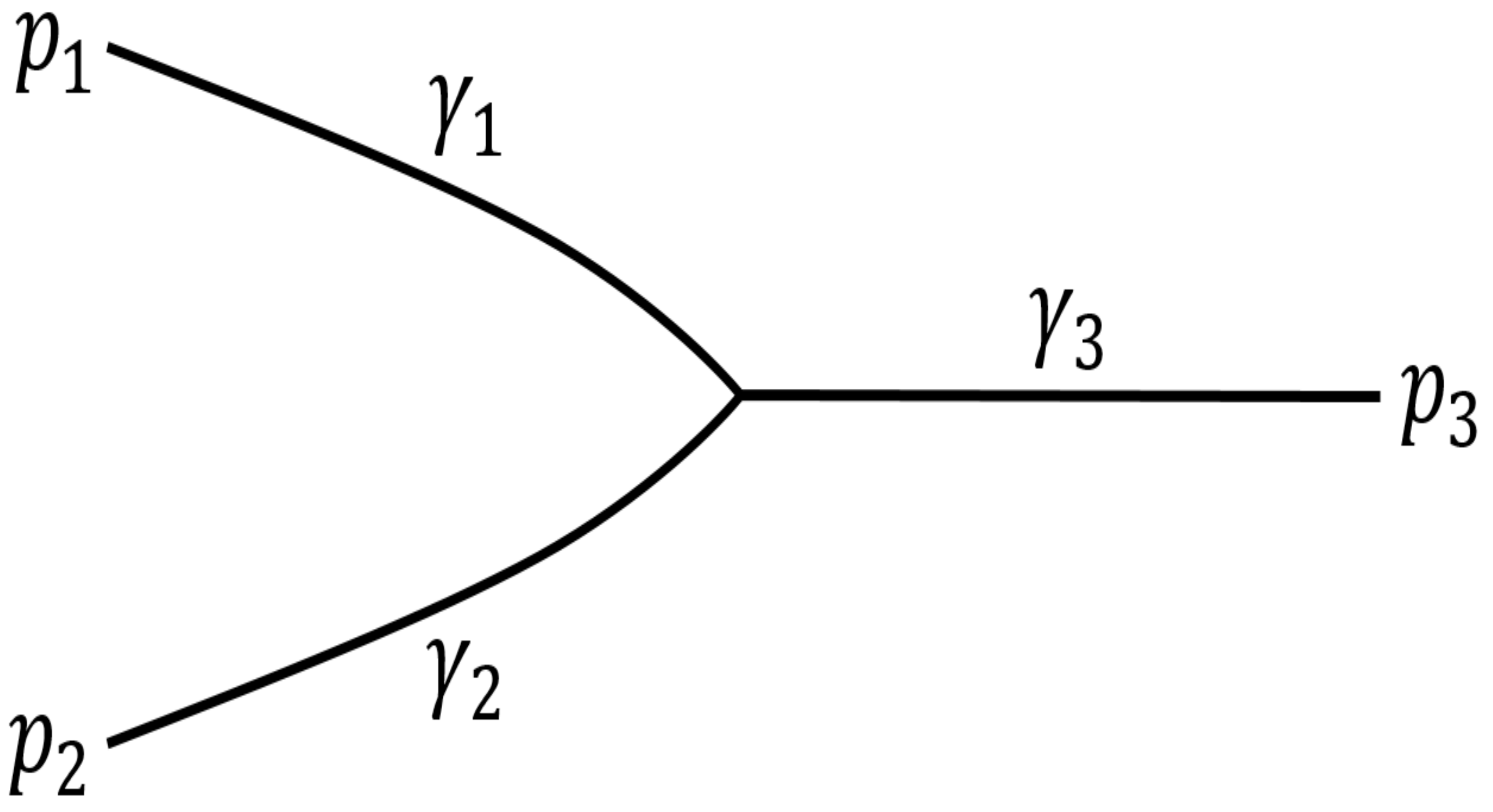}
\centering
\caption{Tree that defines the product $\cdot$}
\label{Tree_pic}
\end{figure}

\noindent{\it Proof.}
All the products are defined using the corresponding products on homology groups for approximations (or for a fixed Riemannian metric). More precisely, we first define:
$$\begin{aligned}
&\circ:HF_*(O_M,\Upsilon:H_1,J_\Upsilon)\otimes HF_*(O_M,\Upsilon:H_2,J_\Upsilon)\longrightarrow HF_*(O_M,\Upsilon:H_3,J_\Upsilon)\\
&\cdot:HM_*(f_N^1:g)\otimes HM_*(f_N^2:g)\longrightarrow HM_*(f_N^3:g)\\
&\star:HM_*(f_N:g)\otimes HF_*(O_M,\Upsilon:H,J_\Upsilon)\longrightarrow HF_*(O_M,\Upsilon:H,J_\Upsilon).
\end{aligned}$$
The above product are defined in a standard way, in the chain level, using the mappings with domains depicted in Figures~\ref{fig:pants},~\ref{Tree_pic} and

\begin{figure}
\centering
\includegraphics[width=8cm,height=4cm]{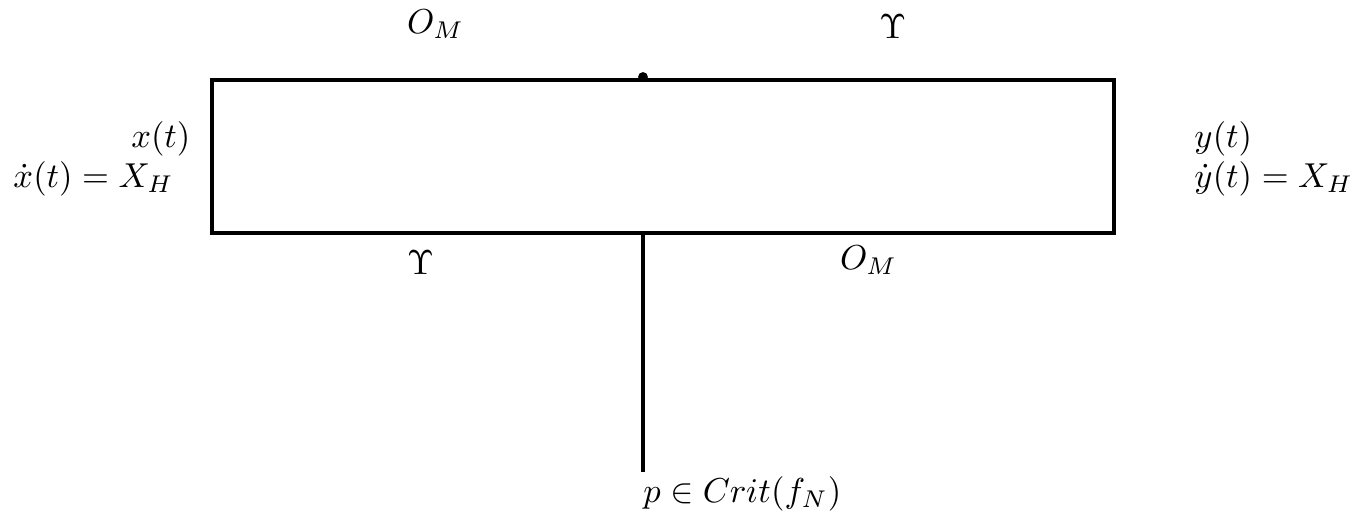}
\centering
\caption{Mixed object that defines the product $\star$}
\label{fig:star}
\end{figure}

Next we check that it holds:
\begin{equation}\label{eq:com+F+prod}
\begin{aligned}
&\mathbf{F}^{H_3}_{ab}(x_a\circ y_a)=\mathbf F^{H_1}_{ab}(x_a)\circ \mathbf F_{ab}^{H_2}(y_a)\\
&\mathbf{G}_{ab}^{f_3}(p_a\cdot q_a)=\mathbf{G}_{ab}^{f_1}(p_a)\cdot\mathbf{G}_{ab}^{f_2}(q_a)\\
&\mathbf{F}_{ab}(p_a\star x_a)=\mathbf{G}_{ab}(p_a)\star\mathbf{F}_{ab}(x_a)
\end{aligned}
\end{equation}
which enables us to have the products well defined on the direct limit levels. The subscript letters $f_j$ and $H_j$ indicate which Morse or Hamiltonian functions we have in mind.

Finally, we check:
\begin{equation}\label{eq:PSS-structure}
\begin{aligned}
&\Phi^{\Upsilon}_3(\alpha\cdot\beta)=\Phi_1^{\Upsilon}(\alpha)\circ\Phi_2^{\Upsilon}(\beta)\\
&(p\cdot q)\star x=p\star(q\star x),
\end{aligned}
\end{equation}
for all $p,q\in HM_*(f_N:g)$ and $x\in HF_*(O_M,\Upsilon:H,J_\Upsilon)$. This gives us the claimed properties of products, on the approximation level, hence, from~(\ref{eq:com+F+prod}), the claim properties on direct limit level.

The details of all above proofs are the same as in~\cite{KMN}, since the constructions apply to the approximations.\qed

\subsection{Triangle inequality}

For two function $H_1, H_2:T^*M\times[0,1]\to\mathbb{R}$ with $H_1(x,1)=H_2(x,0)$, we define their concatenation as:
$$H_1\sharp\,H_2:=\begin{cases}H_1(x,t),&t\le 1\\
H_2(x,t-1),&t\ge 1.
\end{cases}$$

\begin{thrm}\label{thrm:tirangle} Let $f_j$, for $j=1,2,3$ be three admissible Morse functions on $N$. Let $[\alpha]\in HM_*(f_1)$, $[\beta]\in HM_*(f_2)$ be such that $HM_*(f_3)\ni[\alpha\cdot\beta]\neq0$. It holds
$$
c^3_N([\alpha\cdot\beta],H_1\sharp H_2)\le c^1_N([\alpha],H_1)+c^2_N([\beta],H_2),
$$
where $j$ in $c_N^j$ emphasizes the corresponding Morse function $f_j$.
\end{thrm}

\noindent{\it Proof.} One first proves the triangle inequality for approximations:
$$
c^3_\Upsilon(\alpha\cdot\beta,H_1\sharp H_2)\le c^1_\Upsilon(\alpha,H_1)+c^2_\Upsilon(\beta,H_2).
$$ Choose a regular Hamiltonian $H_3$ with
$$\|H_3-H_1\sharp H_2\|_{C^0}<\varepsilon.$$ Denote by $\Sigma$ the Riemannian surface in Figure~\ref{fig:pants}
$$
{\mathbb R}\times[-1,0]\sqcup{\mathbb R}\times[0,1]
$$
with the identification $(s,0^-)\sim(s,0^+)$ for $s\geq0$, and by $(s,t)$ the coordinates defined in three ends, $\Sigma_1\approx(-\infty,0]\times [0,1]$, $\Sigma_2\approx(-\infty,0]\times [0,1]$ and $\Sigma_3\approx([0,\infty)\times [0,1]$.

Let $K:T^*M\times\Sigma\to{\mathbb R}$ be a smooth family of Hamiltonians such that
$$K(\cdot,s,t)=\begin{cases}H_1(\cdot,t+1),&s\le -1,-1\le t\le0\\
        H_2(\cdot,t),&s\le -1,0\le t\le1\\
        \frac{1}{2}H_3(\cdot,\frac{t+1}{2}),&s\ge1\\
        0, &(s,t)\notin\Sigma\setminus(\Sigma_1\cup\Sigma_2\cup\Sigma_3);
        \end{cases}\quad \left\|\frac{\partial K}{\partial s}\right\|\le\begin{cases}\varepsilon,&s\in[-1,1],\\
        0, &\mathrm{elsewhere}.\end{cases}$$
Let $\tilde x=\Phi_1^\Upsilon(\alpha)$ and $\tilde y=\phi^\Upsilon(\beta)$. Since $\alpha\cdot\beta\neq0$ and
$$0\neq\Phi_3^\Upsilon(\alpha\cdot\beta)=\Phi_1^\Upsilon(\alpha)\circ\Phi_2^\Upsilon(\beta)=\tilde x\circ \tilde y,$$
there exist $x\in CF_*^\lambda(O_M,\Upsilon:H_1,J)$, $y\in CF_*^\sigma(O_M,\Upsilon:H_2,J)$ that participate in the formal sums defining $\tilde{x}$ and $\tilde{y}$ respectively, as well as $z\in CF_*(O_M,\Upsilon:H_3,J)$ and $u:\Sigma\to T^*M$ satisfying
$$\left\{\begin{array}{l}
\bar{\partial}_{K,J^{\Upsilon}}(u)=0\\
u(s,-1)\in O_M,\,u(s,1)\in\Upsilon,\,s\in{\mathbb R}\\
u(s,0^-)\in\Upsilon,\,u(s,0^+)\in O_M,\,s\le 0\\
 u_1(-\infty,t)=x(t)\\
 u_2(-\infty,t)=y(t)\\
 u_3(+\infty,t)=z(t)
\end{array}\right.
$$
Since
$$
\begin{aligned}
&\int_{\Sigma}\bigg\|\frac{\partial u}{\partial s}\bigg\|^2\,ds\,dt\ge0\\
&\int_\Sigma u^*\omega=-\int x^*\theta+h_{\Upsilon}(x(1))-\int y^*\theta+h_{\Upsilon}(y(1))+\int z^*\theta-h_{\Upsilon}(z(1)),
\end{aligned}
$$
it follows form properties of a Hamiltonian $K$ and Stoke's formula:
$$
{\mathcal A}^{\Upsilon}_{H_3}(z)\le{\mathcal A}^{\Upsilon}_{H_1}(x)+{\mathcal A}^{\Upsilon}_{H_2}(y)+4\varepsilon.
$$ The proof now follows from Theorem~\ref{thm:Cont_invar} and~(\ref{eq:PSS-structure}).\qed

\subsection{Comparison with spectral invariants for periodic orbits in $T^*M$}

Let us sketch the construction of spectral invariants for periodic orbit Floer homology. Since Floer homology for periodic orbits is not well defined for compactly supported Hamiltonians in $T^*M$, we need to consider Hamiltonians with a support in some fixed cotangent ball bundle. This is used in~\cite{FS} and also in~\cite{MVZ}. Fix $R>0$, $\varepsilon>0$ and a smooth function $h:(-\varepsilon,+\infty)\to\mathbb{R}$ such that:
\begin{itemize}
\item $h(t)=0$ for $t\ge 0$;
\item $h'(t)\ge 0$ for $t\le 0$;
\item $h'$ is small enough so that the flow of $h(\|\mathbf p\|-R)$ does not have non constant periodic orbit of period less or equal to $1$ for $\|\mathbf p\|\in(0,\varepsilon)$.
\end{itemize}
We choose $H_t(\mathbf{q},\mathbf{p})$ to be equal to $h(\|\mathbf{p}\|-R)$ for $\|\mathbf{p}\|\ge R-\varepsilon$.

Let $HF_*(T^*M:H,J)$ and $HM_*(F,T^*M)$ denote Floer homology for periodic orbits in $T^*M$ and Morse homology for the Morse function $F:T^*M\to\mathbb{R}$ respectively. The filtration in Floer homology for periodic orbits is given by the standard action functional
$$a_H(\gamma):=\int\gamma^*\theta-\int_0^1 Hdt$$
which is well defined in the cotangent bundle setting. Denote by $HF_*^{\lambda}(T^*M:H,J)$ the corresponding filtered group and by $\jmath_*^\lambda$ the map induced by the inclusion map. Let $\mathrm{PSS}$ stands for PSS isomorphism for periodic orbits, defined in a way analogous to~\cite{PSS}:
$$\mathrm{PSS}:HM_*(F,T^*M:g)\stackrel{\cong}{\longrightarrow}HF_*(T^*M:H,J)$$ and let $\alpha\in HM_*(f,T^*M:g)$.

Filtered Floer homology groups $HF_*^{\lambda}(T^*M:H,J)$ are homology groups of a chain complex generated by
$$CF_*^{\lambda}(T^*M:H):=\{a\in CF_*(T^*M:H)\mid a_H(a)<\lambda\},$$
where $CF_*(T^*M:H)$ denotes the $\mathbb Z_2-$vector space over the set of periodic Hamiltonian $H-$orbits in $T^*M$.

Let
$$\jmath^\lambda_*:HF_*^{\lambda}(T^*M:H,J)\to HF_*(T^*M:H,J)$$ denote the map induced by the inclusion map. For $\alpha\in HM_*(F,T^*M:g)\setminus\{0\}$ define
$$\rho(\alpha,H):=
\inf\{\lambda\mid\mathrm{PSS}(\alpha)\in\mathrm{Im}(\jmath^\lambda_*)\}.$$

We now choose a Morse function $F:T^*M\to\mathbb{R}$ in a specific way. Let $f_N$ be an admissible Morse function on $N$ and $f_M$ be its extension to $M$ as described on the page~\pageref{page:extension}. We can extend $f_M$ to $F:T^*M\to\mathbb{R}$ in the same way as we extended $f_N$ to $F_M$ to obtain the Morse function $F$ on $T^*M$ with no negative gradient trajectories leaving $N$. Now the Morse complex $CM_*(f_N)$ is a subset of the Morse complex $CM_*(F)$ with $m_{f_N}(p)=m_F(p)$, for $p\in\mathrm{Crit}(f_N)$, and the inclusion map of these complexes becomes the homomorphism $\imath_*$ on the homology level.

For two generic almost complex structures $J_a$ and $J_b$, denote by $\mathbf{D}_{ab}$ a canonical isomorphism of Floer homologies for periodic orbits:
$$\mathbf{D}_{ab}:HF_*(T^*M:H,J_a)\to HF_*(T^*M:H,J_b)$$
that satisfies
$$\mathbf{D}_{bc}\circ \mathbf{D}_{ab}=\mathbf{D}_{ac}.$$ We define a preorder on the space of compatible complex structures to be the full relation and Floer homology for periodic orbits as a direct limit
$$HF_*(T^*M:H):=\displaystyle\lim_{\longrightarrow} HF_*(T^*M:H,J_s).$$ Similarly, by varying a Riemannian metrics, we define Morse homology as a direct limit
$$HM_*(F,T^*M):=\displaystyle\lim_{\longrightarrow} HM_*(F,T^*M:g).$$

It is not hard to check that the inclusion map
$$\imath_*:HM_*(f_N:g)\to HM_*(F,T^*M:g)$$ induces the map, denote it the same:
$$\imath_*:HM_*(f_N)\to HM_*(F,T^*M).$$

\begin{thrm}\label{thm:chimney} Let $[\alpha]\in HM_*(f_N)\setminus\{0\}$. Then
$$c_N([\alpha],H)\ge\rho(\imath_*([\alpha]),H).$$

\end{thrm}

\noindent{\it Proof.} The proof is divided to several steps, we only sketch them here and refer the reader to Section~4 in~\cite{KMN} for the detailed construction and technicalities.

\begin{figure}
\centering
\includegraphics[width=8cm,height=4cm]{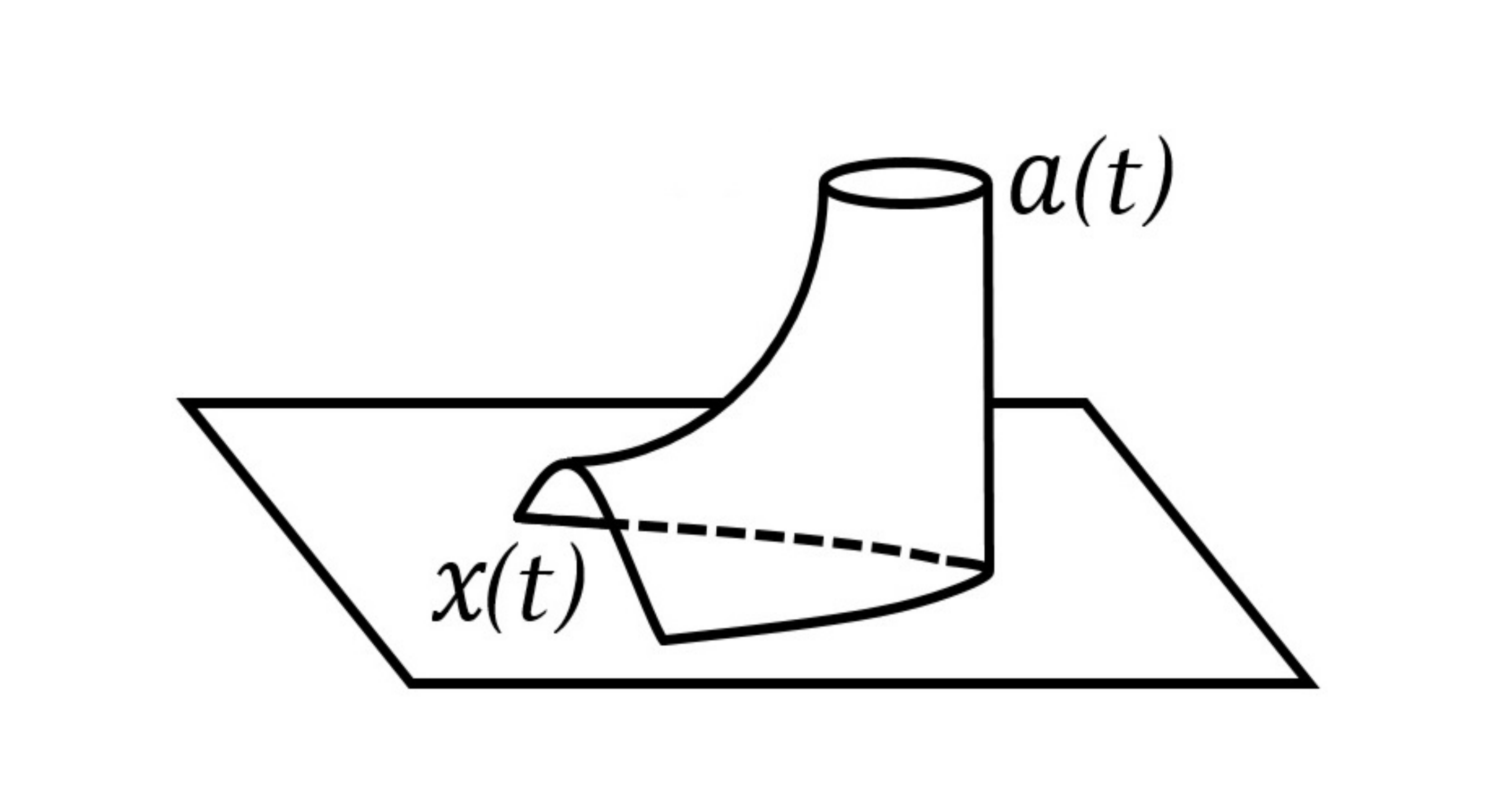}
\centering
\caption{Chimney}
\label{fig:chimney}
\end{figure}

The first step is to construct an Albers-type morphism:
$$\chi:HF_*(O_M,\Upsilon:H,J)\to HF_*(T^*M:H,J),$$ for fixed parameters and approximations. This mapping is defined via the number of perturbed holomorphic ``chimneys", i.e.\ maps $u$ defined on
$$\Pi:=\mathbb{R}\times [0,1]/\sim,\quad\mbox{where}\; (s,0)\sim (s,1)\;\mbox{for}\;s\ge 0$$ that satisfy:
$$
\left\{\begin{array}{l}u:\Pi\to T^*M\\
\partial_su+J(\partial_tu-X_H\circ u)=0\\
u(s,0)\in O_M,\,u(s,1)\in\Upsilon\;\mbox{for}\;s\le 0\\
u(-\infty,t)=x(t), \,u(+\infty,t)=a(t),
\end{array}\right.
$$
see Figure~\ref{fig:chimney} (see also~\cite{A}).

The next step is to show that the maps $\chi$ are also well defined on the filtered groups:
$$\chi^\lambda:HF_*^{\lambda}(O_M,\Upsilon:H,J_\Upsilon)\to HF_*^{\lambda}(T^*M:H,J_\Upsilon).$$ This is easily done by inspection the change of the action functional along chimney.

Then, one shows that the diagrams
\begin{equation}\label{2diags}
\xymatrix{ HF_*^{\lambda}(O_M,\Upsilon:H,J_\Upsilon)\ar[d]_{\imath_*^{\lambda}}\ar[r]^-{\chi^{\lambda}} &HF_*^{\lambda}(T^*M:H,J_\Upsilon)
\ar[d]_{\jmath_*^{\lambda}}\\
HF_*(O_M,\Upsilon:H,J_\Upsilon)\ar[r]^-{\chi} &HF_*(T^*M:H,J_\Upsilon)
\ar[d]_{\mathrm{PSS}^{-1}}\\
HM_*(f_N:g_\Upsilon)\ar[u]_{\Phi^\Upsilon}\ar[r]^{\imath_*} &HM_*(F,T^*M:g_\Upsilon)}
\end{equation}
and
\begin{equation}\label{eq:diag_chi_ab}
\begin{array}{lllllllll}
\cdots&\longrightarrow&HF_*^{\lambda}(\Upsilon_a)&\stackrel{\mathbf F_{ab}}{\longrightarrow}&HF_*^{\lambda}(\Upsilon_b)&\stackrel{\mathbf F_{bc}}{\longrightarrow}&HF_*^{\lambda}(\Upsilon_c)&\longrightarrow&\cdots\\
&&\downarrow\chi^a&&\downarrow\chi^b&&\downarrow\chi^c&&\\
\cdots&\longrightarrow&HF^{\lambda}_*(J_a)&\stackrel{\mathbf D_{ab}}{\longrightarrow}&HF^{\lambda}_*(J_b)&\stackrel{\mathbf D_{bc}}{\longrightarrow}&HF^{\lambda}_*(J_c)&\longrightarrow&\cdots
\end{array}
\end{equation}
commute. This is technically the most demanding step and it can be done by a careful choice of auxiliary one-dimensional manifolds.

The commutativity of the diagram~(\ref{eq:diag_chi_ab}) enables us to define a chimney-type map
$$\chi:HF_*(H,U:M)\to HF_*(T^*M:H).$$
The commutativity of the diagram~(\ref{2diags}) implies the commutativity of
$$\xymatrix{ HF_*^{\lambda}(H,N:M)\ar[d]_{\imath_*^{\lambda}}\ar[r]^-{\chi^{\lambda}} &HF_*^{\lambda}(T^*M:H)
\ar[d]_{\jmath_*^{\lambda}}\\
HF_*(H,N:M)\ar[r]^-{\chi} &HF_*(T^*M:H)
\ar[d]_{\mathrm{PSS}^{-1}}\\
HM_*(f,N)\ar[u]_{\Phi}\ar[r]^{\imath_*} &HM_*(F,T^*M)}$$ which easily implies
$$c_N([\alpha],H)\ge\rho(\imath_*([\alpha]),H).$$ \qed

\section{The case of nested open subsets}\label{section:nested}
Let $U$ be an open subset of $M$, with $\partial U$ a smooth $(n-1)$-dimensional manifold. Our construction applies in this situation ($N=U$, $\dim N=\dim M$), in fact, this case is the subject of~\cite{KMN}.

In~\cite{O3} Oh considered a spectral invariant
$$c_+(H,U):=\inf\{\lambda\in\mathbb{R}\mid \imath^\lambda_*:HF^\lambda_*(H,U:M)\to HF_*(H,U:M)\;\mbox{is surjective}\}.$$

If $U\stackrel{\jmath_{UV}}{\hookrightarrow}V$ are two open subset of $M$ and
$$\jmath_{UV*}:H*(U)\to H_*(V)$$ is surjective, Oh proved that
\begin{equation}\label{eq:O-ineq}
c_+(H,V)\le c_+(H,U).\end{equation}
In~\cite{KMN} we proved a slightly more precise statement, the above inequality for any homology class (with $\mathbb{Z}_2$ coefficients), using the PSS isomorphism for an open subset.

\begin{thrm}\label{thm:c_UV} Let $U{\hookrightarrow}V$ be two open subsets of $M$ and let
$$\jmath_{UV*}:HM_*(f,U)\to HM_*(f,V)$$ (the homomorphism induced by inclusion $\jmath_{UV}:U\hookrightarrow V$) be surjective. Let $c_U([\alpha],H)$ be as in~(\ref{eq:def_inv_c_open}). For $[\alpha]\in HM_k(f,U)\setminus\{0\}$ it holds:
$$c_V(\jmath_{UV*}([\alpha]),H)\le c_U([\alpha],H).$$
\end{thrm}
We generalize~(\ref{eq:O-ineq}) in the following way. Define:\label{c_+}
$$c_+(H,N):=\inf\{\lambda\in\mathbb{R}\mid \imath^\lambda_*:HF^\lambda_*(H,N:M)\to HF_*(H,N:M)\;\mbox{is surjective}\}$$ and for a fixed approximation:
$$c_+(H,\Upsilon):=\inf\{\lambda\in\mathbb{R}\mid \imath^\lambda_*:HF^\lambda_*(O_M,\Upsilon:H,J)\to HF_*(O_M,\Upsilon:H,J)\;\mbox{is surjective}\}.$$

\begin{rem}\label{rem:c+const}
It easily follows from Proposition~\ref{prop:limit} that the invariants $c_+(H,\Upsilon^N)$ all become equal to $c_+(H,N)$ for all approximations $\Upsilon^N\ge\Upsilon_0$.
\end{rem}

The natural arising question is whether spectral invariants are in some way continuous with respect to a submanifold $N\subset M$? We will give an affirmative answer in the case when $N$ is a framed submanifold in $M$. We say that a submanifold $N\subset M$ is {\it framed} if for some tubular neighbourhood $\mathrm{Tb}_N(\partial N)$, the normal bundle $\nu(N\cup\mathrm{Tb}_N(\partial N))$ of a submanifold $N\cup\mathrm{Tb}_N(\partial N)$ is trivial.\label{framed}

\begin{thrm}\label{thrm:cont} Suppose that a submanifold $N\subset M$ is framed.
Let $U_n$ be a decreasing sequence of open subsets with smooth boundaries and $\bigcap_nU_n=N$, where $N$ is a smooth submanifold with a smooth boundary. Then $\lim_{n\to\infty}c_+(H,U_n)=c_+(H,N)$.
\end{thrm}

Before proving Theorem~\ref{thrm:cont}, we first note that we can assume that $U_n$ are tubular neighbourhoods of $N$. Indeed, otherwise, since $U_n$ are open and $N$ is compact, we can find a sequence $V_n$ of tubular neighbourhoods, such that, for $n\ge n_0$ it holds $V_n\subseteq U_n\subseteq V_{n+1}$.

We will divide the proof into several steps.

The first step is to construct, for a fixed $n\in\mathbb{N}$ and $\varepsilon>0$,
an inclusion mapping from $HF_*^\lambda(O_M,\Upsilon^N:H,J)$ to $HF_*^{\lambda+\varepsilon}(O_M,\Upsilon^n:H,J)$, where the approximations $\Upsilon^N$ of $\nu^*\overline{N}$ and $\Upsilon^n$ of $\nu^*\overline{U}_n$ are suitably chosen.

This is done in Lemma~\ref{lem:D_R->D_delta}, Lemma~\ref{lem:rotation}, Lemma~\ref{lem:estimate} and Lemma~\ref{lem:homology_level}. The proofs of Lemma~\ref{lem:D_R->D_delta} and Lemma~\ref{lem:rotation} follow the ideas from~\cite{O3} up to some point, with some specificities adapted to our case. Therefore we expose them in details.

Let $\|\cdot\|$ be the norm induced by a fixed Riemannian  metric on $T^*M$. Denote by
$$D^*_aM:=\{(\mathbf{q},\mathbf{p})\in T^*M\mid\|\mathbf{p}\|\le a\}.$$
Let $R>0$ be such that all the solutions of~(\ref{d}) are contained in $D_R^*M$ (this $R$ exists by standard $C^0$-estimates).

\begin{lem}\label{lem:D_R->D_delta} Fix $\varepsilon>0$, $\delta>0$ and an approximation $\Upsilon^N$ of $\nu^*\overline{N}$. There exists a smooth function $K:T^*M\to\mathbb{R}$ such that the Hamiltonian flow $\phi^1_K$ has the following properties:
\begin{itemize}
\item $\phi^1_K$ leaves $O_M$ and $\nu^*\overline{N}$ invariant
\item $\phi^1_K$ maps $\nu^*\overline{N}\cap D_R^*M$ into $\nu^*\overline{N}\cap D_\delta^*M$
\item $\phi^1_K$ leaves $(T^*M\setminus D^*_{2R}M)\cup(T^*M\setminus T^*(N\cup\mathrm{Tb}_N(\partial N)))$ fixed.
\item $C^0$ norm of $K$ is smaller than $\varepsilon$.
\end{itemize}

\end{lem}

\noindent{\it Proof.} Denote
$$X:=N\cup\mathrm{Tb}_N(\partial N).$$
Since $N$ is framed, so
\begin{equation}\label{eq:trivialN}
\mathrm{Tb}_M(X)\cong X\times\mathbb{R}^{n-k}\;\Rightarrow \;
 T^*\big(\mathrm{Tb}_M(X)\big)\cong T^*(X)\times\mathbb{R}^{n-k}\times\mathbb{R}^{n-k},
\end{equation}
we can choose coordinates
$$(q_1,\ldots,q_n,p_1,\ldots,p_n)\in\mathbb{R}^{2n}$$ that satisfy both~(\ref{eq:coordinates}) and
\begin{itemize}
\item $(q_{k+1},\ldots,q_n,p_{k+1},\ldots,p_n)$ are globally defined
\item $(q_k,p_k)$ are globally defined in $T^*(\mathrm{Tb}_N(\partial N))$, say for $|q_k|<\delta_0$.
\end{itemize}

Recall that, in the coordinates~(\ref{eq:coordinates}), the set $\nu^*_-(\partial N)$ is locally described by
$$\{(q_1,\ldots,q_{k-1},0,\ldots,0,0,\ldots,0,p_k,\ldots,p_n)\mid p_k\le 0\}$$ and $\nu^*N$ is described as
$$\{(q_1,\ldots,q_{k},\underbrace{0,\ldots,0}_{n-k},\underbrace{0,\ldots,0}_{k},p_{k+1},\ldots,p_n)\mid q_k\le 0\}.$$
Let us first define a Hamiltonian diffeomorphism in the set $|q_k|<\delta_0/2$.

Consider a Hamiltonian vectorfield
$$(q_k,\ldots,q_n,p_k,\ldots,p_n)\mapsto(\lambda q_k,\ldots,\lambda q_n,-\lambda p_k,\ldots,-\lambda p_n)$$ in the $(q_k,\ldots,q_n,p_k,\ldots,p_n)$-plane generated by the Hamiltonian function
$$(q_k,\ldots,q_n,p_k,\ldots,p_n)\mapsto\lambda q_kp_k+\ldots+\lambda q_np_n,$$ where $\lambda$ is chosen to be large enough so that the corresponding Hamiltonian diffeomorphism maps the set
$$\{(0,\ldots,0,p_k,\ldots,p_n)\mid \|(0,\ldots,0,p_k,\ldots,p_n)\|\le R\}$$
into
$$\{(0,\ldots,0,p_k,\ldots,p_n)\mid \|(0,\ldots,0,p_k,\ldots,p_n)\|\le\delta\}.$$
Outside $|q_k|<\delta$, our diffeomorphism will be constructed similarly, but without the coordinates $(q_k,p_k)$.

Therefore choose a cut-off function $\beta:[0,+\infty)\to[0,1]$ such that
$$\beta(t)=\begin{cases}0,&t\in[2,+\infty)\\
1,&t\in[0,1].\end{cases}$$
For $a>0$, denote by $\beta_a(t):=\beta\left(\frac{t}{a}\right)$ and define
\begin{equation}\label{eq:K}{K}(\mathbf{q},\mathbf{p}):=\rho_R(\mathbf{q},\mathbf{p})\beta_r(\|(q_{k+1},\ldots,q_n)\|)(\beta_{\delta_0/2}(|q_k|)\lambda q_kp_k+\ldots+\lambda q_np_n).
\end{equation}
Here $\rho_R:T^*M\to[0,1]$ is a smooth cut--off function equal to $1$ in $D^*_RM$ and equal to zero in $T^*M\setminus D^*_{2R}M$
(In~(\ref{eq:K}), as in the rest of the paper, by $\|(q_k,\ldots,q_n)\|$ we assume $\|((0,\ldots,0,q_k,\ldots,q_n,0,\ldots,0)\|$ etc.) By choosing $r>0$ small enough we can obtain $\|K\|_{C^0}<\varepsilon$.
\qed

\medskip

In the following lemma we keep the notations from Lemma~\ref{lem:D_R->D_delta}.

\begin{lem}\label{lem:rotation} For given $\varepsilon>0$, there exist an approximation $\Upsilon^N$ and $\delta>0$ and a Hamiltonian $L$ such that the corresponding Hamiltonian isotopy $\phi_L^1$ satisfies
\begin{itemize}
\item[(i)]  $\phi_L^1$ maps $\Upsilon^N\cap D_\delta^*M$ into $O_M\cap\mathrm{Tb}_M(N\cup\mathrm{Tb}_N(\partial N))\subset O_{U_n}$
\item[(ii)]  $\phi_L^1$ leaves $T^*M\setminus D_{2\delta}^*M$ fixed
\item[(iii)] $\|L\|_{C^0}<\varepsilon$.
\end{itemize}
\end{lem}

\noindent{\it Proof.}  Recall $\Upsilon^N$ is of the form
$$\Upsilon^N:=\sigma(\nu_-^*(\partial N)\cap\nu^*N)\times C$$
(see~(\ref{Upsilon_epsilon})) where $C$ is a smooth curve in $(q_k,p_k)$, Suppose that $\delta>0$ is small enough and $\Upsilon^N$ is such that
$$\Upsilon^0\cap D_\delta\subset\mathrm{Tb}_M(N\cup\mathrm{Tb}_N(\partial N))\quad\mbox{and}\quad C\setminus \left(\{q_k=0\}\cup\{p_k=0\}\right)\subset\mathrm{Tb}_M(N\cup\mathrm{Tb}_N(\partial N)).$$
Since $\Upsilon^N$ decomposes in two components, we will construct a Hamiltonian as a sum of two Hamiltonians:
\begin{align*}
&L(q_1,\ldots,q_n,p_1\ldots,p_n)=L_1(q_k,p_k)+\\
&L_2(q_1,\ldots,q_{k-1},q_{k+1},\ldots,q_n,p_1,\ldots,p_{k-1},p_{k+1},\ldots,p_n).
\end{align*}

First we construct $L_1$. We want $\phi^1_{L_1}$ to
\begin{itemize}
\item rotate the negative part of $p_k$-axis to the positive part of $q_k$-axis
\item leave the line $({q_k}-\mathrm{axis})\cap\Upsilon^N$ fixed.
\end{itemize}
To get this, we can choose the Hamiltonian function of the form $-\frac{\pi}{4}(q_k^2+p_k^2)$ in the region $q_k\le 0$, $p_k\le 0$ and multiply it with a suitable cut-off function, such that it has a compact support and it is equal to zero at $q_k-\mathrm{axis}\cap\Upsilon^N$. Denote the corresponding Hamiltonian diffeomorphism by $\varphi$. It deforms the curve $C$ into some smooth curve, say $\widetilde{C}$ which coincide with $q_k$-axis outside a compact interval. If the curve $\widetilde{C}$ is not a graph of some map $p_k=f(q_k)$, we can deform $C$ in such way that both $C$ and $\widetilde{C}$ become a graph. Note that the Hamiltonian diffeomorphism, say $\psi$ with the differential of corresponding Hamiltonian function equal to $f$ maps $\widetilde{C}$ to the $q_k$-axis. Again we cut-off this Hamiltonian suitably and define $\phi_{L_1}^1$ to be the composition of these two Hamiltonian isomorphisms, $\psi\circ\varphi$.

Let us now construct the Hamiltonian $L_2$. Similarly to the construction above, we first define a Hamiltonian $\check{L}$ in $(q_{k+1},\ldots,q_n,p_{k+1},\ldots,p_n)$-plane such that associated Hamiltonian diffeomorphism which
\begin{itemize}
\item rotates the negative part of $p_j$-axis to the positive part of $q_j$-axis
\item leave the line $({q_j}-\mathrm{axis})\cap\Upsilon^N$ fixed
\end{itemize}
for all $j=k+1,\ldots,n$. This can be done by choosing
$$\check{L}(q_{k+1},\ldots,q_n,p_{k+1},\ldots,p_n):=-\frac{\pi}{4}(q_{k+1}^2+p_{k+1}^2+\ldots+q_n^2+p_n^2)$$
in the region $q_j\le 0$, $p_j\le 0$, $j={k+1},\ldots,n$ and by multiplying it with a suitable cut-off functions in $(q_j,p_j)$-planes, such that it has a compact support and it is equal to zero at $({q_j}-\mathrm{axis})\cap\Upsilon^N$.
Finally we set
\begin{align*}&L_2(q_{k+1},\ldots,q_n,p_{k+1},\ldots,p_n):=\\
&\beta_\delta\left(\|(p_1,\ldots,p_{k-1},p_{k+1},\ldots,p_n)\|\right)\check{L}(q_{k+1},\ldots,q_n,p_{k+1},\ldots,p_n).
\end{align*}
By decreasing $\delta>0$ we obtain (iii).
\qed

\begin{lem}\label{lem:estimate} Denote
\begin{equation}\label{eq:comp_K}
\begin{aligned}
&\tilde{x}(t):=\phi^t_L\circ\phi^t_K(x(t))\\
&\tilde{u}(s,t):=\phi^t_L\circ\phi^t_K(u(s,t))\\
&\widetilde{H}:=L\;\sharp\;K\sharp\;H\\
&\widetilde{J}_t:=\left(\phi^t_L\circ\phi^t_K\right)_*J_t.
\end{aligned}
\end{equation}
Then, for $\delta>0$ small enough, the mapping $x\mapsto\tilde{x}$ maps the set $CF(O_M,\Upsilon^N:H,J)$ into $CF(O_M,\Upsilon^{n}:\widetilde{H},\widetilde{J})$, where $\Upsilon^n$ is an approximation of $\nu^*\overline{U}_n$. The mapping $u\mapsto\tilde{u}$ maps the set $\mathcal{M}(O_M,\Upsilon^N:H,J)$ into $\mathcal{M}(O_M,\Upsilon^n:\widetilde{H},\widetilde{J})$. Moreover, for given $\varepsilon>0$, we can decrease $\delta$ and choose approximations $\Upsilon^n$ and $\Upsilon^N$ such that for every $x\in\mathcal{M}(O_M,\Upsilon^N:H,J)$ it holds
\begin{equation}\label{eq:estimateA_H}
\left|\mathcal{A}^{\Upsilon^n}_{\tilde{H}}(\tilde{x})-\mathcal{A}^{\Upsilon^N}_{H}(x)\right|<\varepsilon.
\end{equation}
\end{lem}

\noindent{\it Proof.} For a fixed $\Upsilon^N$ and $\Upsilon^n$, one can choose $\delta$ small enough such that
$$x\in CF(O_M,\Upsilon^N:H,J)\;\Rightarrow\;\tilde{x}\in CF(O_M,\Upsilon^n:\widetilde{H},\widetilde{J}).$$ This is possible since the set $CF(O_M,\Upsilon^N:H,J)$ is finite.
Then, one easily checks
$$u\in\mathcal{M}(O_M,\Upsilon^N:H,J)\;\Rightarrow\;\tilde{u}\in\mathcal{M}(O_M,\Upsilon^n:\widetilde{H},\widetilde{J}).$$

To prove~(\ref{eq:estimateA_H}), we estimate
\begin{equation}\label{eq:triangle}\left|\mathcal{A}^{\Upsilon^n}_{\tilde{H}}(\tilde{x})-\mathcal{A}^{\Upsilon^N}_{H}(x)\right|\le\left|\mathcal{A}_{\tilde{H}}(\tilde{x})-\mathcal{A}_{H}(x)\right|+|h_{\Upsilon^n}(\tilde{x}(1))|+|h_{\Upsilon^N}({x}(1))|,
\end{equation} where
$$\mathcal{A}_H(\gamma):=\int \gamma^*\theta-\int_0^1H(\gamma(t),t)dt.$$
The last two terms in~(\ref{eq:triangle}) can be made smaller than $\varepsilon/4$ by choosing good enough approximations.

Next, as in~\cite{O3}, note that
$$\mathcal{A}_{\tilde{H}}(\tilde{x})-\mathcal{A}_{H}(x)=\int_0^1\frac{d}{ds}\left[\mathcal{A}_{s(L\,\sharp\,K)\,\sharp\,H}\left(\phi^t_{s(L\,\sharp\,K)}(x)\right)\right]$$
and
$$\frac{d}{ds}\left[\mathcal{A}_{s(L\,\sharp\,K)\,\sharp\,H}\left(\phi^t_{s(L\,\sharp\,K)}(x)\right)\right]=d\mathcal{A}_{s(L\,\sharp\,K)\,\sharp\,H}\left(\frac{\partial\phi^t_{s(L\,\sharp\,K)}(x)}{\partial s}\right)-\int_0^1L\,\sharp\,K(\phi^t_{s(L\,\sharp\,K)}(x))dt.$$
By decreasing $\delta$ and $r$ from~(\ref{eq:K}), we can obtain
$$\left|\int_0^1L\,\sharp\,K(\phi^t_{s(L\,\sharp\,K)}(x))dt\right|<\frac{\varepsilon}{4}.$$

To finish the proof, recall
$$d\mathcal{A}_H(\gamma)(\xi)=\int_0^1\omega(\dot{\gamma}-X_H(\gamma),\xi)dt+\langle\xi,\theta(\gamma(1))\rangle-\langle\xi,\theta(\gamma(0))\rangle.$$ Since
$$s\mapsto\phi^t_{s(L\,\sharp\,K)}(x(t))$$ is a Hamiltonian orbit of $L\,\sharp\,K$ with the initial point $x(t)$ and
$$\left\langle\frac{\partial\phi^0_{s(L\,\sharp\,K)}(x(0))}{\partial s},\theta\left(\phi^0_{s(L\,\sharp\,K)}(x(0))\right)\right\rangle=0$$
we have
\begin{equation}\label{eq:dA_H}
d\mathcal{A}_{s(L\,\sharp\,K)\,\sharp\,H}\left(\frac{\partial\phi^t_{s(L\,\sharp\,K)}(x)}{\partial s}\right)=\left\langle\frac{\partial\phi^1_{s(L\,\sharp\,K)}(x(1))}{\partial s},\theta\left(\phi^1_{s(L\,\sharp\,K)}(x(1))\right)\right\rangle.
\end{equation}
Let us check that the right side of~(\ref{eq:dA_H}) can be made arbitrarily small, by decreasing $\delta$.
Since
$$\frac{\partial\phi^1_{s(L\,\sharp\,K)}(x(1))}{\partial s}=X_{L\,\sharp\,K}$$ which is basically $X_K+X_L$, let us estimate the form $\theta$ at the latter vector fields. Firstly, since the support of $L$ is contained in $D^*_{2\delta}M$, by decreasing $\delta$, we can obtain
$$|\mathbf{p}d\mathbf{q}(X_L)|<\frac{\varepsilon}{4}.$$
Secondly, we have
$$\begin{aligned}
&\frac{\partial K}{\partial p_j}=\beta_r(\|(q_{k+1},\ldots,q_n)\|)\left(\frac{\partial\rho_R}{\partial p_j}(\mathbf{q},\mathbf{p})(\beta_{\delta_0/2}(|q_k|)\lambda q_kp_k+\ldots+\lambda q_np_n)+
\rho_R(\mathbf{q},\mathbf{p})\cdot a(\mathbf{q},\mathbf{p})\right)\\
&\mbox{where}\quad a(\mathbf{q},\mathbf{p})):=\begin{cases}\beta_{\delta_0/2}(|q_k|)\lambda q_k,&j=k
\\
\lambda q_j,&j=k+1,\ldots,n.
\end{cases}
\end{aligned}$$
 Since
$$\left\|\frac{\partial K}{\partial p_j}\right\|\le\|(q_k,\ldots,q_n)\|\cdot const,$$ we can decrease $r$ so we obtain
$$|\theta(X_K)|=|\mathbf{p}d\mathbf{q}(X_K)|=\left|\sum p_j\frac{\partial K}{\partial p_j}\right|<\frac{\varepsilon}{4}.$$\qed

\begin{lem}\label{lem:homology_level} It holds
\begin{equation}
\label{eq:comm_d}\widetilde{\partial x}=\partial\tilde{x}.
\end{equation} and, for every $\varepsilon>0$, the mapping
\begin{equation}\label{eq:i_n}
\imath_n^\lambda=\imath^\lambda_{(\Upsilon^N,\Upsilon^n)}:HF^\lambda(O_M,\Upsilon^N:H,J)\to HF^{\lambda+\varepsilon}(O_M,\Upsilon^n:H,J)
\end{equation}
 is well defined.
\end{lem}

\noindent{\it Proof.} The proof of~(\ref{eq:comm_d}) is similar to the proof in Chapter 3 in~\cite{O3}. We sketch the key steps for the sake of completeness.

For, $\delta>0$, denote by
$$U_{\delta}:=\{\mathbf{q}\in\mathrm{Tb}_M(N\cup\mathrm{Tb}_N(\partial N))\mid\|(q_{k+1},\ldots,q_n)\|<\delta\},$$ where $(q_{k+1},\ldots,q_n)$ are global coordinates from~(\ref{eq:trivialN}).
For a exact Lagrangian manifold $\Upsilon$, denote by
$$\mathcal{M}(O_M,\Upsilon:H,J):=
\left\{ u\left|
\begin{array}{l}
u:\mathbb{R}\times[0,1]\rightarrow T^*M\\
\frac{\partial u}{\partial s}+
J(\frac{\partial u}{\partial t}-X_H(u))= 0 \\
u(s,0)\in O_M,\,u(s,1)\in{\Upsilon}.
\end{array}\right.\right\}.$$
There exists $\delta>0$ and $\lambda$ form Lemma~\ref{lem:D_R->D_delta} such that every $\tilde u\in\mathcal{M}(O_M,\Upsilon^n:\widetilde{H},\widetilde{J})$ satisfies:
\begin{equation}\label{eq:U_delta}
\tilde u(-\infty,1)\in U_{\delta}\Rightarrow \tilde u(s,1)\in U_{\delta},\quad\mbox{for all}\;s\in\mathbb{R}.
\end{equation}
The proof of this claim is the same as the proof of Proposition 3.2 in~\cite{O3}. Note that at some point the proof of  Proposition 3.2 in~\cite{O3} uses the fact that $H\equiv 0$ near $t=0$ and $t=1$. This is not a loss of generality, since, according to Lemma 3.1 in~\cite{O3}, we can perturb our Hamiltonian to satisfy this condition with the change in the filtration as small as desired. The condition~(\ref{eq:U_delta}) implies~(\ref{eq:comm_d}). Indeed choose $\lambda$ such that
$$\tilde{u}(s,1)\in U_\delta,\quad\mbox{for all}\;u\in\mathcal{M}(O_M,\Upsilon^N:H,J).$$ Together with~(\ref{eq:U_delta}), this means that for every $u\in\mathcal{M}(O_M,\Upsilon^N:H,J)$, $\tilde{u}$ participates in the computation of $d\tilde{x}$. Again from~(\ref{eq:U_delta}) and the fact the mappings
$$x\mapsto\tilde{x},\quad u\mapsto\tilde{u}$$ are bijections from $\mathcal{M}(O_M,\Upsilon^N:H,J)$ to the set
$$\{u\in\mathcal{M}(O_M,\Upsilon^n:\widetilde{H},\widetilde{J})\mid u(s,1)\in U_{\delta}\} $$
(since it is constructed using Hamiltonian diffeomorphism), we conclude the reverse: every $\tilde{u}$ participating in $d\tilde{x}$ comes from some $u\in\mathcal{M}(O_M,\Upsilon^N:H,J)$, for every $\tilde{x}$ with $\tilde{x}(1)\in U_\delta$. This implies~(\ref{eq:comm_d}).

Now we compose our map $x\mapsto\tilde{x}$ with a canonical continuation homomorphism (see, for example~\cite{KO}):
$$HF_*^\lambda(O_M,\Upsilon^n:\widetilde{H},\widetilde{J})\to HF^{\lambda+\varepsilon/2}_*(O_M,\Upsilon^n:{H},{J})$$ and finish the proof. This canonical isomorphism exists for a given $\varepsilon$ since the norms $\|K\|_{C^0}$ and $\|L\|_{C^0}$ can be made arbitrary small.\qed

\begin{lem}\label{lem:inequality} For $U_n$ a tubular neighbourhood of $N$, it holds:
$c_+(H,N)\ge c_+(H,U_n)$
\end{lem}
\noindent{\it Proof.} The following diagram
$$\xymatrix{  HF^\lambda_*(O_M,\Upsilon^N:H,J)\ar[d]_{\jmath_{N*}^\lambda}\ar[r]^-{\imath^\lambda_n} & HF_*^{\lambda+\varepsilon}(O_M,\Upsilon^n:H,J)
\ar[d]_{\jmath_{U_n*}^\lambda}\\
HF_*(O_M,\Upsilon^N:H,J)\ar[d]_{F_{(H,N)}^{-1}}\ar[r]^-{\imath_n} &HF_*(O_M,\Upsilon^n:H,J)\ar[d]_{F_{(H,U)}^{-1}}\\
H_*(N)\ar[r]^{\jmath_{NU_n*}}&H_*(U)}
$$
commutes. Indeed, the upper diagram is obvious, and the commutativity od the lower one is proved in~\cite{O3} in the similar situation (for $U\subset V$ open). The vertical mappings $F_{(H,N)}$ and $F_{(H,U)}$ are isomorphisms between singular and Floer homologies defined in~\cite{KO} (without the use of PSS).

Since $N$ is a deformation retract of $U_n$, the map
$$\jmath_{NU_n*}:H_*(N)\to H_*(U_n)$$ induced by the inclusion $\jmath_{NU_n}:N\to U_n$ is an isomorphism. Note that in this case, the mapping $\imath_n$ for $\lambda=\infty$ in~(\ref{eq:i_n}) is also an isomorphism.

Now we easily compute:
$$c_+(H,\Upsilon^n)\le c_+\left(H,\Upsilon^N\right)+\varepsilon.$$

Since the invariants $c_+$ all become equal, starting from some approximation (Remark~\ref{rem:c+const}), we conclude that
$$c_+(H,U)\le c_+(H,N)+\varepsilon$$ for every $\varepsilon$, so the proof is complete.

Note that in this lemma we only used the surjectivity of $\jmath_{NU_n*}$, and not the injectivity. \qed

\bigskip

\noindent{\it Proof of Theorem~\ref{thrm:cont}.} In the same way as in the proof of Lemma~\ref{lem:inequality} we can prove that $c_+(H,U_n)\le c_+(H,U_{n+1}).$ Therefore, we need to show
that, given $\varepsilon>0$ there exists $n_0$ such that
$$c_+(H,N)<c_+(H,U_{n_0})+\varepsilon.$$

Choose $n_0$ such that $N$ is a deformation retract of $U_{n_0}$ and
$\delta>0$, $\Upsilon^N$, $\Upsilon^{n_0}$ such
that, for $\tilde{H}$ and $\tilde{x}$ from Lemma~\ref{lem:estimate} it holds
\begin{itemize}
\item $\left|\mathcal{A}^{\Upsilon^{n_0}}_{\tilde{H}}(\tilde{x})-\mathcal{A}^{\Upsilon^N}_{H}(x)\right|<\frac{\varepsilon}2$ (Lemma~\ref{lem:estimate})
\item$ c_+\left(H,\Upsilon^N\right)=c_+(H,N),\quad c_+\left(H,\Upsilon^{n_0}\right)=c_+(H,U_{n_0})$ (Remark~\ref{rem:c+const})
\item $c_+(\tilde{H},\Upsilon^{n_0})<c_+({H},\Upsilon^{n_0})+\frac{\varepsilon}{2}$.
\end{itemize}
The last item can be achieved by choosing $\delta$ and $r$ from~(\ref{eq:K}) small enough so that $\|H-\tilde{H}\|_{C^0}$ norm is small enough; this fact can be proved by tracking the change of the action functional associated to the homotopy of Hamiltonians (see Theorem 7.2 in~\cite{O1} or Proposition 4.4 in~\cite{KO}).

Note that
$$c_+(H,\Upsilon^N)=\max\{\mathcal{A}_H^{\Upsilon^N}(x)\mid [x]\in HF_*(O_M,\Upsilon^N:H,J)\}$$
and similarly for $c_+(H,\Upsilon^{n_0})$. For every $[x]\in HF_*(O_M,\Upsilon^N:H,J)$, we have
$$\mathcal{A}_H^{\Upsilon^N}(x)<\mathcal{A}^{\Upsilon^n}_{\tilde{H}}(\tilde{x})+\frac{\varepsilon}{2}\le c_+(\tilde{H},\Upsilon^{n_0})+\frac{\varepsilon}{2}\le c_+({H},\Upsilon^{n_0})+\varepsilon.$$ By taking a maximum over  $[x]\in HF_*(O_M,\Upsilon^N:H,J)$ we conclude
$$c_+\left(H,\Upsilon^N\right)<c_+\left(H,\Upsilon^{n_0}\right)+\varepsilon.$$

\qed

\begin{rem}
In the same way we can prove that the invariants\label{c_-}
 $$c_-(H,N):=\sup\{\lambda\in\mathbb{R}\mid \imath^\lambda_*:HF^\lambda_*(H,N:M)\to HF_*(H,N:M)\;\mathrm{is\; trivial}\}$$
are continuous with respect to a decreasing sequence $U_n$, $\bigcap U_n=N$, therefore the same holds for  the invariants:
$$\gamma(H,N):=c_+(H,N)-c_-(H,N).$$
\end{rem}

\begin{ques} An interesting question is whether there is a similar continuity result for spectral invariants for a single homology class, defined in~(\ref{eq:def_inv_c_open}). If this is true, then the inequality
$$c_U(\imath^M_*[\alpha],H)\le c_N([\alpha],H)$$ would easily follow from the commutativity of the following diagram:
$$\begin{array}{lllll}
HM_*(f_N:g)&\stackrel{\mathrm{PSS}^N}{\longrightarrow}&HF_*(O_M,\Upsilon^N:H,J)&\stackrel{\jmath^N_\lambda}{\longleftarrow}&HF_*^{\lambda}(O_M,\Upsilon^N:H,J)\\
\downarrow\imath_*^M&&\downarrow\imath_{(\Upsilon^N,\Upsilon^U)}&&\downarrow\imath_{(\Upsilon^N,\Upsilon^U)}^\lambda\\
HM_*(f_{U}:g)&\stackrel{\mathrm{PSS}^U}{\longrightarrow}&HF_*(O_M,\Upsilon^U:H,J)&\stackrel{\jmath^U_{\lambda+\varepsilon}}{\longleftarrow}&HF_*^{\lambda+\varepsilon}(O_M,\Upsilon^U:H,J).
\end{array}$$
\end{ques}

\begin{rem}\label{rem} The previous results apply to the case when the boundary of a submanifold is empty. Consider two closed framed submanifolds (with or without boundary) $N_1$ and $N_2$ such that $N_1\subset N_2$.
Choose two sequences of open sets $U_n$ and $V_n$ that satisfy:
\begin{itemize}
\item $U_n$ is a tubular neighbourhood of $N_2$, $V_n$ is a tubular neighbourhood of $N_1$
\item $U_{n+1}\subseteq U_n$, $V_{n+1}\subseteq V_n$
\item $\bigcap_nU_n=N_2$, $\bigcap_nV_n=N_1$
\item $\partial U_n$ and $\partial V_n$ are smooth submanifolds of $M$ of codimension one
\item $V_n\subset U_n$.
\end{itemize}
Suppose that
$$\imath_*:H_*(N_1)\to H_*(N_2)$$ is surjective, then $\imath_*:H_*(V_n)\to H_*(U_n)$ is surjective.
It follows from Theorem~\ref{thm:c_UV} that
$$c_+(H,V_n)\le c_+(H,U_n).$$ By taking the limit when $n\to\infty$ and using Theorem~\ref{thm:Cont_invar} we conclude
$$
c_+(H,N_1)\le c_+(H,N_2).
$$
If the inclusion map $\imath_*:H_*(N_1)\to H_*(N_2)$ is injective (hence $\imath_*:H_*(V_n)\to H_*(U_n)$ is injective), Oh showed that
$$
c_-(H,V_n)\ge c_-(H,U_n)
$$
(see Theorem 4.3 (2) in~\cite{O3}). Therefore
\begin{equation}\label{eq:ineq-c-}
c_-(H,N_1)\ge c_-(H,N_2).
\end{equation}
Finally, if $\imath_*:H_*(N_1)\to H_*(N_2)$ is a bijection, it holds
$$\gamma(H,N_1)\le\gamma(H,N_2).$$

Consider the special case when $N_2=N$, $N_1=\partial N$. Suppose that $N$ is framed, this implies that $\partial N$ is framed too, since
$$\nu(\partial N)\cong\nu(N)|_{\partial N}\oplus\mathrm{Tb}_N(\partial N).$$ If $\imath_*:H_*(\partial N)\to H_*(N)$ is surjective, then:
$$c_+(H,\partial N)\le c_+(H,N).$$

Since we deal with $\mathbb{Z}_2$-coefficients, then $\imath_*:H_*(\partial N)\to H_*(N)$ is never an injection, for $*=k=\dim N$ (when $N$ and $\partial N$ are connected). Indeed, since
$$H_k(N)=\{0\},\quad H_k(N,\partial N)=\mathbb{Z}_2,$$ using the long exact sequence for the pair $(N,\partial N)$ we have that
$$\mathrm{Ker}(\imath_*)=\mathbb{Z}_2,\quad \imath_*:H_{k-1}(\partial N)\to H_{k-1}(N).$$
However, if $N$ is not orientable and if we deal with $\mathbb{Z}$-coefficients, the above reasoning does not apply, and it is possible for $\imath_*$ to be injective. In order to prove the inequality~(\ref{eq:ineq-c-})
one needs to consider a coherent orientation of mixed moduli spaces and construct a PSS morphism from Theorem~\ref{thrmPSS-intro} with $\mathbb{Z}$-coefficients (see also~\cite{Frol,KM1}). This will be the subject of further research.
\end{rem}


\end{document}